\documentclass{elsarticle}

\usepackage{lineno,hyperref}
\modulolinenumbers[5]

\usepackage{amsmath,amsfonts,amssymb}
\usepackage{color}
\usepackage{psfrag}
\usepackage{indentfirst}
\usepackage{enumerate}
\usepackage{placeins}
\usepackage{multirow}
\usepackage{booktabs}
\usepackage{array}
\usepackage{mathabx}
\usepackage{algorithm}
\usepackage{algpseudocode}
\usepackage{natbib}
\usepackage{subcaption}
\usepackage{makecell}

\allowdisplaybreaks[4]

\def\v#1{\mbox{\boldmath $#1$}}        % vettore in grassetto
\def\vh#1{\hat {\v {#1}}}

\def\Re{{\mathop{\Bbb R}}}

\makeatletter
\def\BState{\State\hskip-\ALG@thistlm}
\makeatother

\makeatletter
\def\ps@pprintTitle{%
   \let\@oddhead\@empty
   \let\@evenhead\@empty
   \let\@oddfoot\@empty
   \let\@evenfoot\@oddfoot
}
\makeatother
\begin{document}

\begin{frontmatter}

\title{Parameter Estimation of Fire Propagation Models\\ Using Level Set Methods}

\tnotetext[mytitlenote]{This work was supported by the Air Force Office of Scientific Research under Grant FA9550-15-1-0530.}

\author[DIME]{Angelo Alessandri}
\ead[url]{alessandri@dime.unige.it}
\author[DIME]{Patrizia Bagnerini}
\ead[url]{bagnerini@dime.unige.it}
\author[CNR]{Mauro Gaggero\corref{cor1}}
\ead[url]{mauro.gaggero@cnr.it}
\author[DIME]{Luca Mantelli}
\ead[url]{luca.mantelli@edu.unige.it}
\address[DIME]{Department of Mechanical Engineering, University of Genoa, I-16145 Genoa, Italy}
\address[CNR]{Institute of Marine Engineering, National Research Council of Italy, I-16149 Genoa, Italy}
\cortext[cor1]{Corresponding author.}

\begin{abstract}
The availability of wildland fire propagation models with parameters estimated in an accurate way starting from measurements of fire fronts is crucial to predict the evolution of fire and allocate resources for firefighting. Thus, we propose an approach to estimate the parameters of a wildland fire propagation model combining an empirical rate of spread and level set methods to describe the evolution of the fire front over time and space. The estimation of parameters in the rate of spread is performed by using fire front shapes measured at different time instants as well as wind velocity and direction, landscape elevation, and vegetation distribution. Parameter estimation is performed by solving an optimization problem, where the objective function to be minimized is the symmetric difference between predicted and measured fronts at different time instants. Numerical results coming from by the application of the proposed method are reported in two simulated scenarios and in a case study based on data originated by the 2002 Troy fire in Southern California. The obtained results showcase the effectiveness of the proposed approach both from qualitative and quantitative viewpoints.
\end{abstract}


\end{frontmatter}

\section{Introduction}
\label{sec:intro}

Wildland fires are a serious threat for human beings and protected natural areas. Hence, the capability to forecast the evolution of fires over time and space and evaluate the hazard level of a certain site is crucial. Many different modeling approaches have been developed for this purpose, including full-physics models, which provide the most accurate results by taking into account many of the phenomena involved in fire propagation together with their mutual interactions \cite{Hanson2000TheSimulation}. Such models consider regional and local weather as well as terrain characteristics to forecast fire behavior, and also predict production and dispersion of smoke. However, full-physics models may be too computationally demanding for real-time applications, short-term forecasting, and parameter estimation. As an alternative, empirical models have been developed in the literature, such as the one proposed by Rothermel \cite{Rothermel1972AFuels}, where the rate of spread of the fire is described through correlations based on experimental results. The complexity of empirical models may be very different, ranging from models where wind is constant \cite{BurganASubsystem} to models that account
for local interactions between heat sources and wind \cite{Clark1996ADynamics}. Indeed, the level of detail adopted to describe topological features and vegetation (or fuel) distribution may vary so much over space
to undermine simulation accuracy \cite{LIDAR}.

Based on the knowledge of the rate of spread, the evolution of the fire front can be described by means of different techniques, which can be split into two main categories
depending on the implementation, i.e., raster and
vector. Examples of raster techniques are those based on cellular automata \cite{Karafyllidis1997AAutomata, French1990GraphicalSpread, Ghisu2015AnSpread,
HernandezEncinas2007ModellingAutomata, CellularAutomata08Greek, tonini2020}, while the most popular vector approaches are level set methods \cite{Munoz-Esparza2018AnMethod, Mentrelli2016ModellingMethod, Mallet2009ModelingMethods}. In raster implementations, the fire evolution only depends on the interaction among contiguous cells, which can be denoted as burnt, burning, or not burning. The fire propagation  from one cell to its neighbors is defined by a set of rules. Such methods are computationally efficient, but the fire front can be subject to significant distortions. Many solutions have been proposed to mitigate this issue, such as the increase of the number of possible spread directions \cite{French1990GraphicalSpread}, particular kinds of discretization of the spatial domain (e.g., using hexagonal grids) \cite{HernandezEncinas2007ModellingAutomata}, and the adoption of suitable correction factors \cite{Ghisu2015AnSpread}.
In vector implementations, the fire front is defined explicitly by a given amount of points or implicitly by the solution of a partial differential equation (PDE). In both cases, the fire shape is updated at different time steps based on the rate of spread. These methods are accurate, but are more computationally demanding than raster implementations. This issue worsens in the presence of merging fire fronts or unburnt areas when using explicit methods (sometimes referred to as Lagrangian), but it is less significant for implicit methods (also called Eulerian) such as level set methods. The reader interested in a deeper comparison between Lagrangian and Eulerian methods in the context of fire propagation is referred to \cite{Bova} and \cite{Kaur}. Owing to their properties, level set methods are widely employed also in many other application fields, like image processing, computational fluid dynamics, and material science \cite{sethian, ss03, fedkiw, kimmel}. Instead, their use in the context of estimation and control has received much less attention from the research community (among the available results, see \cite{BernHerz11, YangTomlCDC13, AleBagGagTNNLS19, AleBagCiaGag2019, AleBagGagRosAUT20}).

In this paper, we focus on a wildland fire propagation model based on the combination of an empirical description of the rate of spread and the use of level set methods to account for fire front evolution. For this model, we propose an approach for the optimal estimation of the main parameters that results from the minimization of a least-squares cost fitting the available measures of the fire front at different time instants. In our case, the cost to be minimized is not continuously differentiable, and therefore we have to resort to non-derivative methods such as the generalized pattern search (GPS) algorithm \cite{Bertsekas2016NonlinearProgramming}. After obtaining a fire propagation model with accurately estimated parameters, the evolution of the fire front can be predicted in real time, thus allowing the use of the model in a decision support system for fire extinguishing.

The estimation of parameters involved in wildland fire propagation models starting from available measures is a topic that has attracted the interest of researchers in the last decades. To this purpose, various terminology has been used, such as ``model calibration'' or ``parameter identification.'' As pointed out in  \cite{Lautenberger13}, wildland fire model calibration is typically accomplished manually, and only very recently optimization algorithms have been adopted to automate the estimation process, either for Lagrangian \cite{LagrangianParameterEstimation} or Eulerian \cite{Lautenberger13} models. In more detail, \cite{ParametricUncertainty} focuses on parametric uncertainty quantification of the rate of spread. Moreover, references are available that perform parameter estimation using Kalman filtering approaches \cite{Kalman1, Kalman2, Kalman3}. Other studies combine parameter and state estimation based on a Luenberger observer framework \cite{StateParameterEstimation} or adopt other data-driven techniques based on fire propagation simulators, also exploiting historical data \cite{DataAssimilation, DataDriven}.
Such methods perform dynamic estimation over time and are quite complex. Further, they are based on specific fire propagation models that are different from the one considered in this work. Instead, the method presented in this paper proposes a simple yet effective method for calibration of the parameters that relies only on observations of the fire front over time and wind measurements to perform estimation of several parameters (some of them depending on vegetation) that are involved in the rate of spread. Toward this end, the minimization of a cost function measuring the symmetric difference between the observed and predicted fronts is performed\footnote{This idea was exploited by the authors also in \cite{AleBagGagTNNLS19, AleBagCiaGag2019, AleBagGagRosAUT20} with the goal of finding optimal control policies for propagating fronts (in different settings than wildland fires). Here the context is completely different since we focus on parameter estimation rather than control.}. The  considered fire propagation model combines the rate of spread proposed in \cite{Lo2012AMethod, Mallet2009ModelingMethods} and level set methods to model front evolution. To the best of the authors' knowledge, no previous attempts concerning parameter estimation for the model adopted in this paper are available in the literature. This model was used also in \cite{Zhai2020} to predict fire evolution in a real-scale shrubland fire scenario using machine learning techniques, but no parameter estimation was performed. However, the focus on the rate of spread reported in \cite{Lo2012AMethod, Mallet2009ModelingMethods} is not a limitation for the  proposed approach, which could be extended also to more complex expressions with no conceptual difficulties. The effectiveness of the proposed method is tested from qualitative and quantitative viewpoints in different scenarios using both synthetic and real data. Another approach focused on parameter estimation based on level set methods is proposed in \cite{Lautenberger13}, but it relies on different equations and parameters as compared to our work. In more detail, the main contribution of \cite{Lautenberger13} is the development of ELMFire, a geospatial model for simulating wildland fires using level set methods. Several equations are inserted into the model to account for surface fire acceleration, crown fire initiation, dead fuel moisture content, spot fire formation, as well as adjustment factors due to wind. Then, a simple genetic optimization approach is presented to automatically calibrate baseline model parameters, by adopting a fitness function based on a scoring policy that compares burned and not burned areas in a pixel-by-pixel fashion, differently from our approach where we evaluate the symmetric difference between simulated and measured front. However, unlike our work, in \cite{Lautenberger13} no quantitative evaluation of the effectiveness of proposed approach is given by comparing several intermediate fronts during fire evolution.

The rest of this paper is structured as follows. Section \ref{sec:model} describes the considered  model for fire propagation. The proposed technique for parameter estimation is presented in Section \ref{sec:ident}. Section \ref{sec:sim} reports numerical results on model testing and parameter estimation in two simulated case studies and in a real one. Conclusions are drawn in Section \ref{sec:concl}.

\section{Wildland Fire Front Propagation Modeling}
\label{sec:model}

In this section, we investigate the propagation of a wildland fire with different environmental conditions by combining an empirical approach to compute the rate of spread and level set methods to represent the moving fire front.

\subsection{Fire Rate of Spread}
\label{sec:spread}

\begin{figure*}[tb]
\centering
\includegraphics[width=6.0cm]{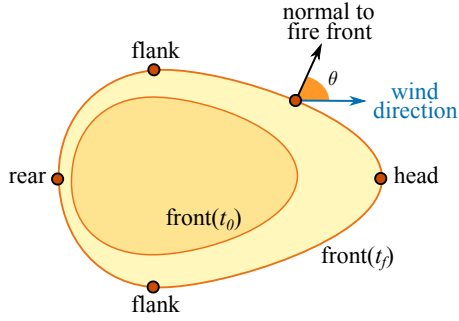}
\caption{Propagation of a wildland fire front from time $t_0$ to $t_f$ under the influence of wind.}
\label{FrontRear.ps}
\end{figure*}

We consider the model of wildland fire rate of spread proposed by Lo in \cite{Lo2012AMethod}, which is derived by previous works of Rothermel \cite{Rothermel1972AFuels} and others \cite{Mallet2009ModelingMethods}.
The main assumption of such model is that, in the presence of wind, if we consider a flat landscape with homogeneous fuel, the fire spread velocity is different at the head, rear, and flanks of the front. As showcased in Fig. \ref{FrontRear.ps}, by ``head'' we mean the region where fire propagation is in the wind direction, by ``rear'' we denote the region where it is in the opposite direction of wind, and by ``flanks'' we indicate the region where it is across of wind. 

\begin{figure*}[tb]
\centering
\subcaptionbox{}{\includegraphics[width=3.8cm]{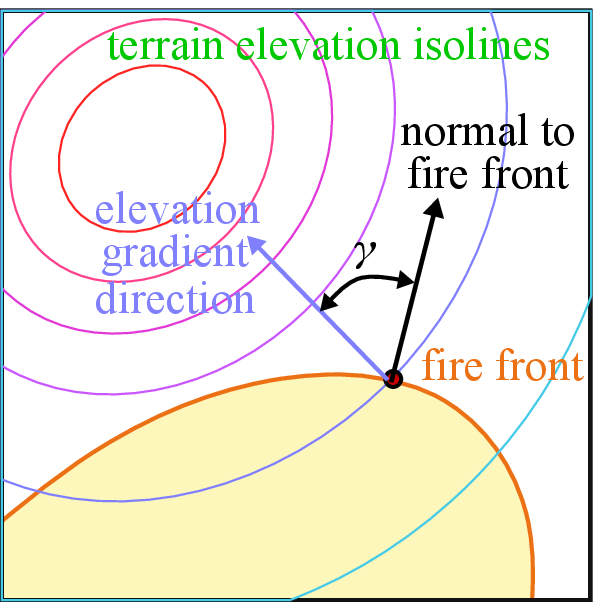}}%
%\hfill
\hspace{0.5cm}
\subcaptionbox{}{\includegraphics[width=3.8cm]{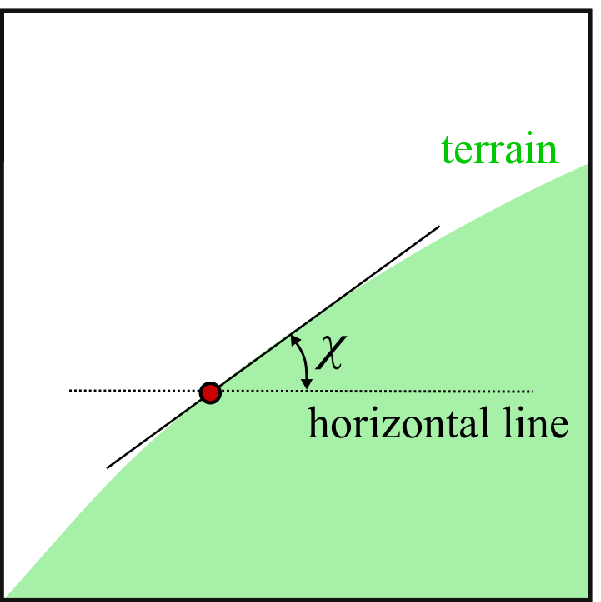}}%
\caption{Sketches of the definitions of the angles $\gamma$ (a) and $\chi$ of the rate of spread (b). A red color of the terrain elevation isolines in (a) denotes a greater elevation with respect to a blue one.}
\label{slope.ps}
\end{figure*}

The head is the area with the highest velocity since wind blows high-temperature combustion products towards unburnt regions. Instead, at the rear such products are driven over already-burnt vegetation, and therefore velocity is lower. At the flanks, the component of the wind speed that is perpendicular to the front is zero, but from experimental tests it has been observed that front propagation is faster than in the absence of wind. However, according to \cite{Lo2012AMethod}, the effect of wind on the spread rate at the flanks is negligible since its impact on the front is small\footnote{The interested reader can find a deeper discussion on the modeling of flank fire using level set methods in \cite{Hilton1, Mentrelli2016ModellingMethod, Hilton2}. A detailed investigation of this aspect is outside the scope of this paper.}. If we assume that a wildland fire propagates only in the normal direction to the front, the front rate of spread $F(U,\theta,\gamma,\chi)$ is given by
\begin{equation}
F(U,\theta,\gamma,\chi) := \begin{cases} 
\!\epsilon + a\sqrt{U \cos^n (\theta)} + \epsilon\psi(\chi) \cos(\gamma), \quad \mbox{if}\:\: \theta\in\Big[0,\frac{\pi}{2}\Big], \\
\epsilon (\alpha+(1-\alpha)|\sin(\theta)|) + \epsilon\psi(\chi) 
\cos(\gamma), \quad \mbox{if}\:\: \theta\in\Big(\frac{\pi}{2},\pi\Big], \end{cases}
\label{eq:F}
\end{equation}
with
\begin{equation}
\psi(\chi) := 5.275\,\beta^{-0.3}\,(\tan(\chi))^2,
\label{eq:Phis}
\end{equation}
where $U\ge 0$ is the module of the wind speed, $\theta\in[0,\pi]$ is the angle between the wind and the normal direction to the front head (see Fig. \ref{FrontRear.ps}), and $\gamma\in[0,\pi]$ is the angle between the elevation gradient direction and the normal to the fire front (see Fig. \ref{slope.ps}(a)). The model also depends on parameters $a\ge 0$ taking into account fuel, $\epsilon\ge 0$  accounting for the velocity of fire at the flanks, $\alpha\ge 0$ representing the ratio between the velocity of fire at the rear and $\epsilon$, and $n\ge 0$ influencing the way the front propagates. The term $\psi(\chi)$ takes into account the terrain effect and includes the fuel bed packing ratio parameter $\beta\ge 0$ and the landscape slope $\tan(\chi)$, where $\chi\in [-\pi/2,+\pi/2]$ is the angle between the terrain and the horizontal line (see Fig. \ref{slope.ps}(b)). Since fire naturally tends to move upward, velocity is higher when propagation is uphill, while it is lower when it spreads downhill \cite{Mallet2009ModelingMethods,Lo2012AMethod}. Toward this end, the term $\psi(\chi)$ in \eqref{eq:F} is multiplied by the cosine of the angle $\gamma$. The coefficients $5.275$ and $-0.3$ in \eqref{eq:Phis} have been fixed according to the reference literature \cite{Rothermel1972AFuels, Lo2012AMethod}.

The accuracy of the model \eqref{eq:F}-\eqref{eq:Phis} strictly depends on the chosen values for the parameters $n$, $\epsilon$, $a$, $\alpha$, and $\beta$. Thus, in Section \ref{sec:ident} we will propose an approach for their optimal estimation based on observations of the fire front. Instead, the quantities $U$, $\theta$, $\gamma$, and $\chi$ in \eqref{eq:F}, \eqref{eq:Phis} will be supposed to be known since they depend on wind and terrain elevation.

\subsection{Level Set Methods}

\begin{figure*}[tb]
\centering
\includegraphics[width=8.0cm]{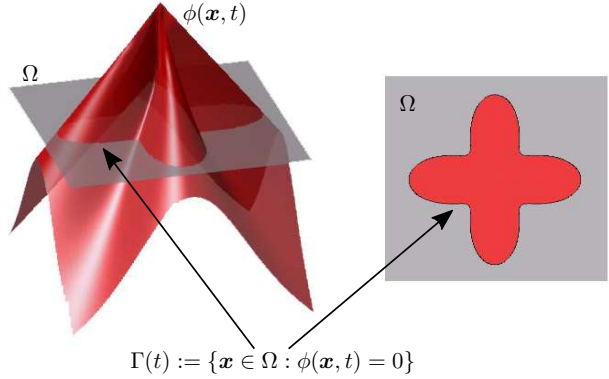}
\caption{Example of function $\phi$ and its zero level set.}
\label{phi.ps}
\end{figure*}

We account for the evolution of the fire front by means of an Eulerian approach 
based on level set methods. Toward this end, let $\Omega \subset \Re^2$ and $t \in [t_0,t_f]$ be a compact two-dimensional space domain and the time,
respectively, where $t_0\ge 0$ and $t_f>t_0$ are given initial and final instants, respectively. Level set methods represent a moving front at each time $t$, i.e., a curve in two 
dimensions separating two regions, as
the zero level set of a function $\phi:\Omega \times [t_0,t_f] \rightarrow
\mathbb{R}$ (see Fig. \ref{phi.ps}). The front $\v x(t,s)$
is given at time $t$ by the points such that $\phi(\v x(t,s),t)=0$,
where $s$ is the arc-length parameter of the initial curve $\v x(t_0,s)$.  
If we differentiate with respect to $t$, we
obtain 
\begin{equation}
\phi_t(\v x,t) + v(\v x,t) \cdot \nabla\phi(\v x,t)=0.
\label{eq:HJeq}
\end{equation}
Equation \eqref{eq:HJeq} is a Hamilton-Jacobi PDE, where $v(\v x,t):=\frac{d}{dt}\v x(t,s)$ is the Lagrangian particle velocity giving the direction of propagation of the front at the point $\v x(t,s)$ and $\nabla$ denotes the gradient with respect to space. In the case of wildland fires, we focus on the normal flow equation since we assume that propagation only occurs in the normal direction to the fire front (see Section \ref{sec:spread}). This corresponds to choose $v(\v x,t)$ proportional to the normal $\v n$ to the front, i.e.,
\begin{equation}
v(\v x,t) := S\, \v n = S \, \frac{\nabla \phi(\v x,t)}{|\nabla \phi(\v x,t)|}, \,
\label{eq:v}
\end{equation}
where $S$ is the speed of front propagation and $n:=\nabla\phi/|\nabla\phi|$ is the normal to the front. In our case, $S$ is given by the rate of spread $F(U,\theta,\gamma,\chi)$ defined in \eqref{eq:F}. If we replace \eqref{eq:v} in \eqref{eq:HJeq} and $S$ with $F(U,\theta,\gamma,\chi)$, we get the following normal flow equation: 
\begin{equation}
\phi_t(\v x,t) + F(U,\theta,\gamma,\chi) \, |\nabla \phi(\v x,t)| = 0.
\label{eq:levelset_Futheta}
\end{equation}
Since $F(U,\theta,\gamma,\chi)$ depends on the angle $\theta$ between the front (i.e., the zero level set of $\phi$) and the wind, it can be considered as a function of $\phi$. Thus, with a little abuse of notation, we can re-write \eqref{eq:levelset_Futheta} as follows:
\begin{equation}
\phi_t(\v x,t) + F(\phi) \, |\nabla \phi(\v x,t)| = 0.
\label{eq:levelset_F}
\end{equation}
Equation \eqref{eq:levelset_F} is associated with initial conditions $\phi_0: \Omega
\rightarrow \Re$, i.e., $\phi(\v x,t_0) = \phi_0(\v x)$ for all $\v x \in \Omega$. 
Usually,
$\phi_0$ is the signed distance to the initial front, and it is conventionally positive inside the fire front (the burnt area) and negative outside (the unburnt area).

At each time $t$, the fire front is represented by the zero level set of the function $\phi$ that solves \eqref{eq:levelset_F}, i.e., it
is a set-valued mapping $\Gamma : [t_0, t_f] \rightrightarrows \mathcal{C}\subset\Re^2$,
where $\Gamma(t) := \left\{
\v x \in \Omega : \phi(\v x,t) =0 \right\}$. 

\begin{figure*}[tb]
\centering
\includegraphics[width=6.0cm]{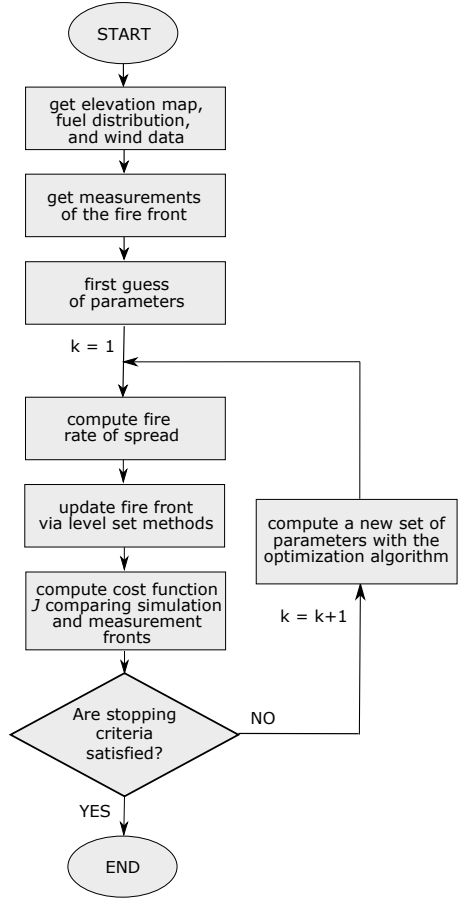}
\caption{Schematic diagram of the parameter estimation procedure.}
\label{fig:block}
\end{figure*}

Since the rate of spread $F$ defined in \eqref{eq:F} depends on the angle between the fire front and the wind, it is not defined everywhere in the domain $\Omega$, but only for the points $\v x$ such that $\phi(\v x,t)=0$. However, it is possible to extend the expression of $F$ in all $\Omega$ by considering the angle formed by all level sets of $\phi$ with the wind. In this way, $F$ is defined for all points $\v x\in\Omega$. The value of $\theta$ in \eqref{eq:F} is computed at each time $t$ as the angle between the normal to the front and the wind direction $\v U$, i.e.,
\begin{displaymath}
\theta = \arccos\left(\frac{\nabla\phi(\v x,t)\cdot \v U}{|\nabla\phi(\v x,t)|\, U}\right),
\end{displaymath} 
where $\cdot$ denotes the usual scalar product and $U$ is the module of the vector $\v U$. 

Level set methods are well-suited to describing wildland fire evolution since no front reconstruction algorithms are required to determine the front shape at a given time. In fact, the fire front can be obtained implicitly from $\phi$, which is defined continuously on the spatial domain at any time instant. For this reason, level set methods are able to easily manage changes of topology, such as merging of fire fronts, or large deformations.

\section{Parameter Estimation}
\label{sec:ident}

The model of fire propagation described in Section \ref{sec:model} depends on several parameters describing the characteristics of the region where fire propagates, such as the type of vegetation or terrain elevation. In this section, we investigate how to compute optimal estimates of  quantities related to the kind of fuel, while the quantities depending on terrain elevation will be assumed known. In more detail, we focus on estimation of the parameters $n$, $\epsilon$, $a$, $\alpha$, and $\beta$ involved in the rate of spread. A precise knowledge of the values of such parameters is crucial to accurately forecast the evolution of wildland fires and therefore plan effective firefighting procedures. 

The estimation of parameters in the rate of spread equation \eqref{eq:F} is
performed by using measurements of the fire front shape at different time instants, velocity and direction of wind, landscape elevation, and fuel distribution. Without loss of generality, the wind is assumed to be constant (as regards both velocity and direction) in the domain $\Omega$ over the whole simulation time. Moreover, we assume that two different kinds of fuels are present in $\Omega$, denoted as ``fuel A'' and ``fuel B.'' We suppose to know the fuel distribution over $\Omega$, i.e., whether fuel A or fuel B is present in the various regions of the domain (such information can be obtained, for instance, from aerial observations \cite{aerialFire}), but its properties are supposed to be unknown. In other words, the values of the parameters $\epsilon$, $a$, and $\beta$ for the two fuels, denoted by $\epsilon^A$, $a^A$, and $\beta^A$ for fuel A and $\epsilon^B$, $a^B$, and $\beta^B$ for fuel B, are assumed to be unknown. The parameters $n$ and $\alpha$ are unknown as well, while the quantities $\theta$, $\gamma$, and $\tan(\chi)$ can be easily computed by combining the available information on the considered region (i.e., landscape elevation) and observations of the fire front. Hence, there is no need to estimate them. Summarizing, we collect all unknown parameters in a $8$-dimensional vector $\v p$, as follows:
\begin{displaymath}
\v p := \left( n, \epsilon^A, \epsilon^B, a^A, a^B, \alpha, \beta^A, \beta^B \right) \in \Re^8.
\end{displaymath}

Let $\Gamma^{\rm meas}(t)$ be the fire front measured from observations at time $t$ and let $\phi^{\rm meas}(x,t)$ be the corresponding level set function obtained as the signed distance to the front (conventionally, positive inside the front and negative outside). Moreover, let $\Gamma(t,\v p)$ be the fire front predicted by the model described in Section \ref{sec:model} corresponding to a given value for the parameter vector $\v p$, and let $\phi(x,t,\v p)$ be the corresponding level set function obtained as the signed distance to the front using the same convention described above. The best estimate $\vh p$ of $\v p$ can be obtained by minimizing a least-squares fitting cost function as follows:
\begin{displaymath}
\vh p := \arg \: \min_{\v p\in [\v p_{\rm min},\v p_{\rm max}]} J(\v p)
\end{displaymath}
where $\v p_{\rm min}\in\Re^8$ and $\v p_{\rm max}\in\Re^8$ are given lower and upper bounds for the vector $\v p$, respectively, and
\begin{equation}
J(\v p) := \int_{t_0}^{t_f} \eta \left( \Gamma(t,\v p) \: \,\Delta \: \,\Gamma^{\rm meas}(t)
\right) \, dt.
\label{eq:Delta}
\end{equation}
The operator $\Delta$ in \eqref{eq:Delta} is the symmetric difference, i.e., $\Gamma(t,\v p) \:\: \Delta \:\: \Gamma^{\rm meas}(t) := (\Gamma(t,\v p) \cup \Gamma^{\rm meas}(t)) \setminus (\Gamma(t,\v p) \cap \Gamma^{\rm meas}(t))$, and $\eta$ corresponds to an outer measure on $\Re^2$. 

\begin{figure*}[tb]
\centering
\subcaptionbox{}{\includegraphics[width=3.1cm]{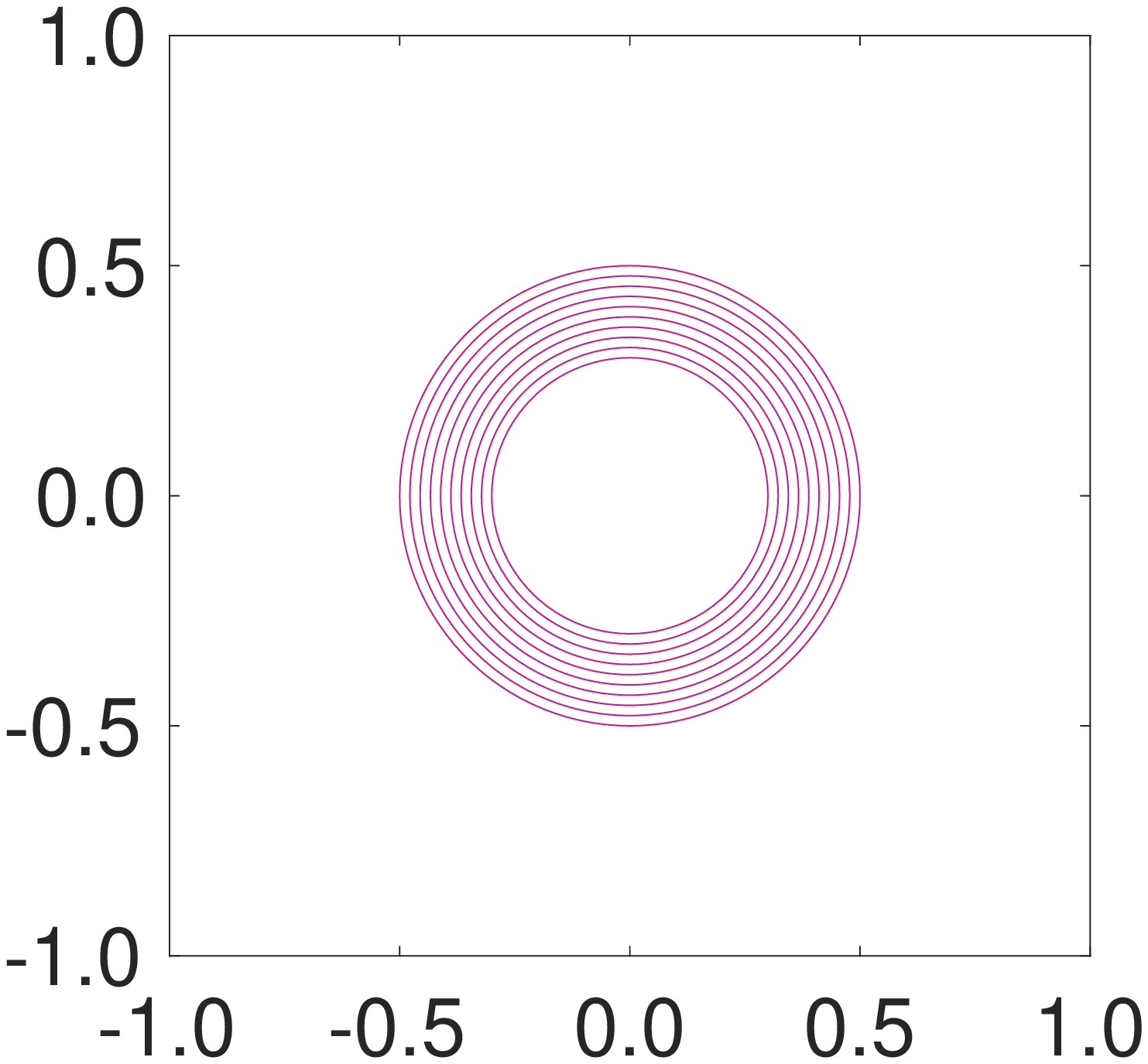}}\label{Sim_NoVento.ps}%
%\hfill
\subcaptionbox{}{\includegraphics[width=3.1cm]{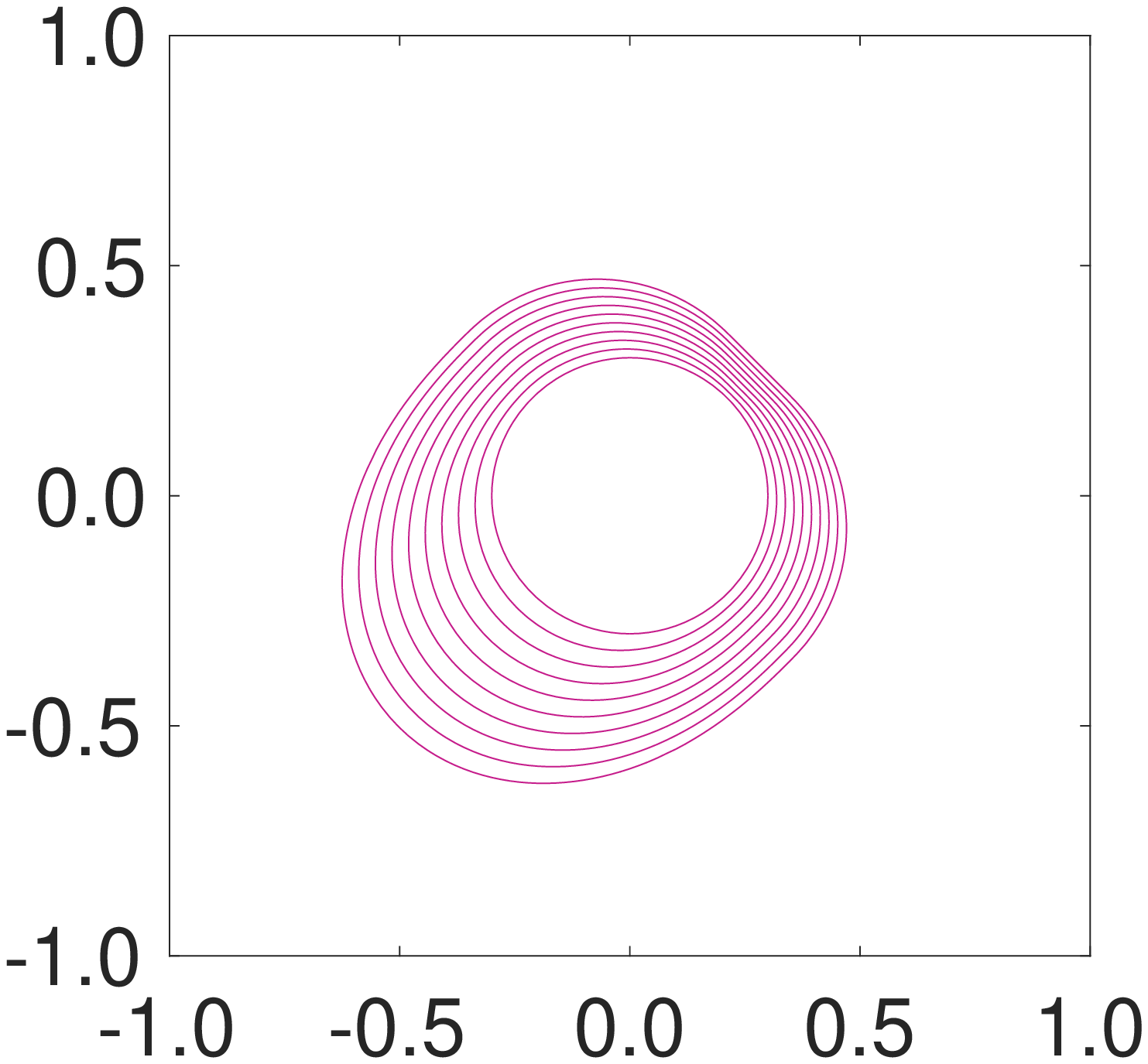}}\label{Sim_Vento.ps}%
%\hfill
\subcaptionbox{}{\includegraphics[width=3.1cm]{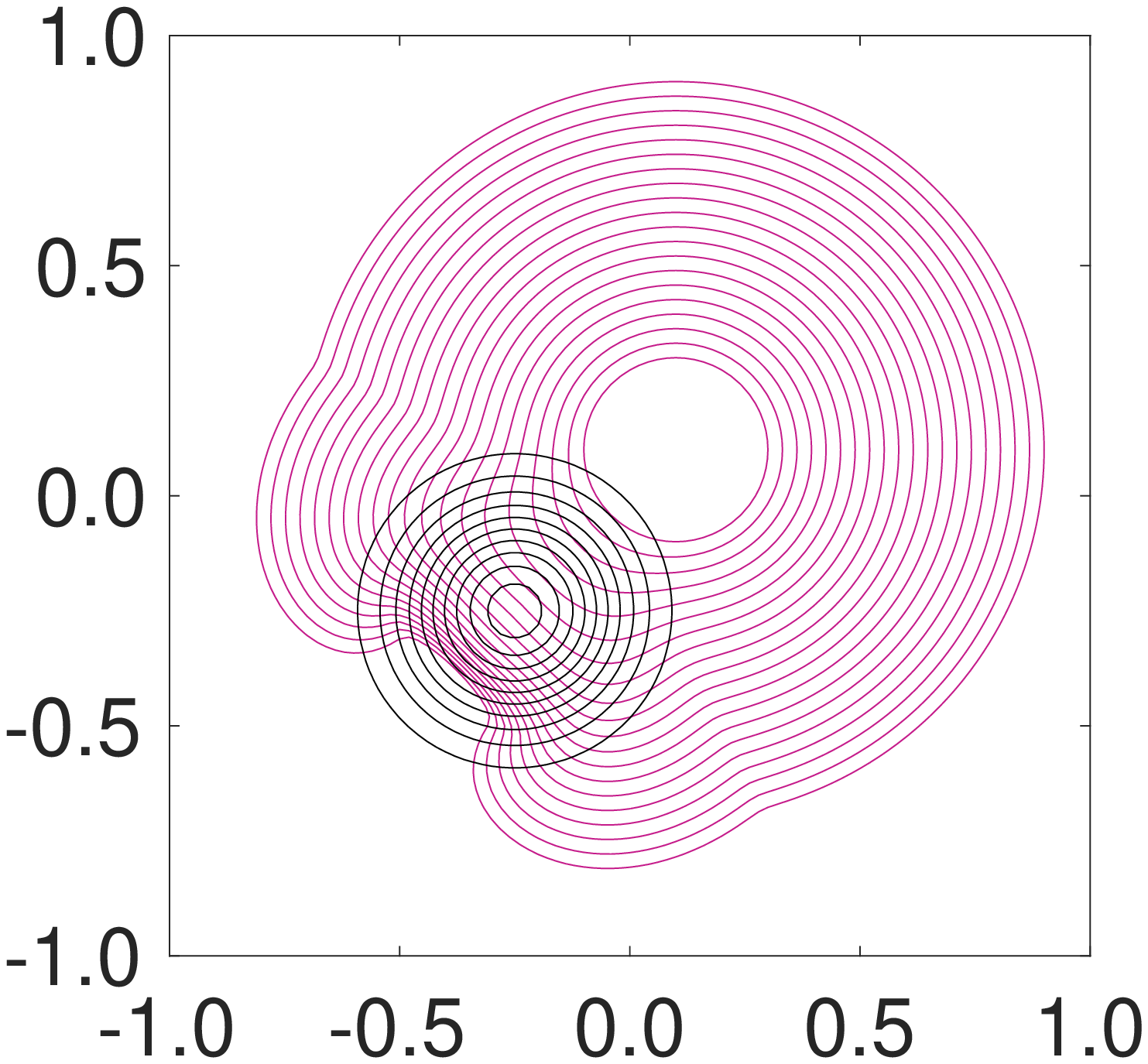}}\label{Sim_Collina.ps}%
%\hfill
\subcaptionbox{}{\includegraphics[width=3.1cm]{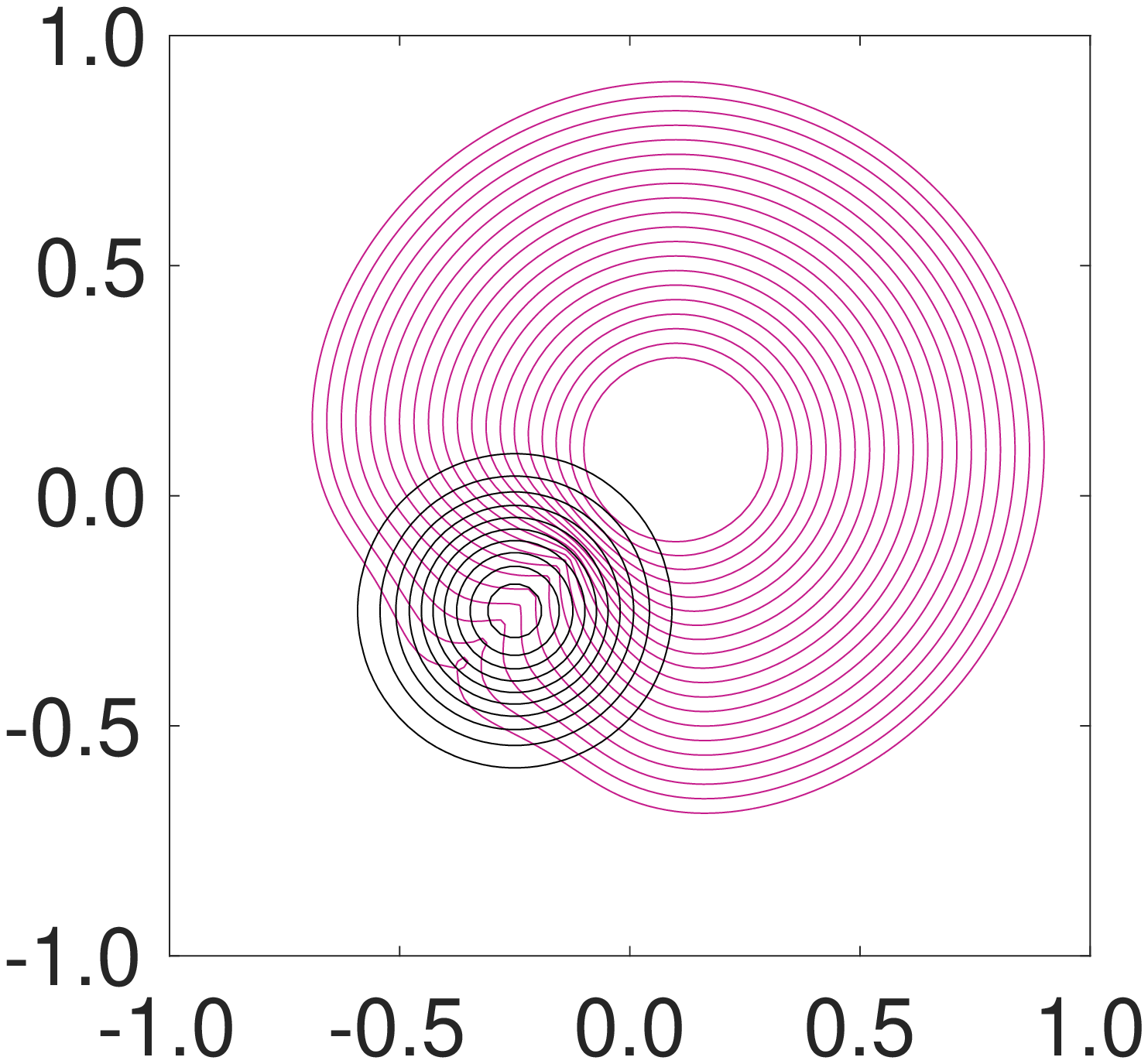}}\label{Sim_Buco.ps}%
\caption{Fire front propagation without wind on a flat landscape (a), with wind on a flat landscape (b), without wind in the presence of a hill (c), and without wind in the presence of a depression (d) (fire fronts are depicted in magenta and elevation isolines are represented in black).}
\label{fig:figsim}
\end{figure*}

Owing to the definition of the symmetric difference and the convention used for the sign of the function $\phi$, the cost  \eqref{eq:Delta} can be rewritten as follows to highlight its dependence on $\phi$:
\begin{equation}
J(\v p) := \int_{t_0}^{t_f} \!\!\int_\Omega \left( H(\phi(\v x,t,\v p)) \!-\!
H(\phi^{\rm meas}(\v x,t)) \right)^2 d\v x \, dt,
\label{eq:trackingCostH}
\end{equation}
where $H(\cdot)$ is the Heaviside step function. 

At each iteration of the optimization algorithm chosen to minimize the cost $J$ in \eqref{eq:trackingCostH}, a certain value for the parameter vector $\v p$ is used to run a fire propagation simulation in the time interval $[t_0,t_f]$ and the domain $\Omega$. Hence, the cost $J$ is computed by comparing the predicted fronts and the measured ones. Then, an updated value for the parameter vector $\v p$ is chosen by the optimization algorithm, and a new simulation of fire propagation is performed to compute a new value for the cost. The procedure is iterated until convergence of the optimization algorithm is obtained and stopping criteria are satisfied. In more detail, they are satisfied if a ``small'' norm of the difference between the estimated parameters or the values of the cost function in two consecutive iterations is found, as well as if a maximum number of iterations is exceeded. Fig. \ref{fig:block} sketches a schematic diagram of the considered parameter estimation procedure, where $k$ denotes a generic iteration of the optimization algorithm used to minimize \eqref{eq:trackingCostH}.

Unfortunately, the cost function associated with the front evolution may be affected by several local minima. Thus, optimization may be subject to local minima trapping, which can provide poor fitting and inaccurate parameter estimates.

\section{Numerical Results}
\label{sec:sim}

In this section, we first evaluate the effectiveness of the model considered in Section \ref{sec:model} to predict fire fronts in simple case studies taken from the reference literature \cite{Lautenberger13}. Then, we focus on parameter identification using the approach presented in Section \ref{sec:ident}. In particular, to evaluate the effectiveness of the proposed technique for estimating the parameter vector $\v p$, we consider three case studies denoted by ``Valley,'' ``Hill,'' and ``Troy.'' The first and second ones are based on simulated data, while the third case involves measurements from a real fire event.

Either for testing or parameter estimation, a suitable discretization method is required for the numerical solution of the model of fire front propagation described in Section \ref{sec:model}. Toward this end, we discretize the spatial domain on a squared regular grid, with cells characterized by given fuel properties and elevation. We compute slope and elevation gradients along the North-South and East-West directions from the elevation matrix using a finite-differences central scheme for inner points and single-sided finite differences for the boundary. 

To numerically solve the Hamilton-Jacobi equation \eqref{eq:levelset_F}, we adopt the Matlab toolbox of level set methods \cite{Mitchell2008TheMethods}. In more detail, an upwind third-order essentially non-oscillatory scheme \cite{kimmel} is applied to approximate space derivatives, while a three-steps second-order Runge-Kutta scheme is used to compute time derivatives. From the  values of $\phi$, it is straightforward to obtain the fire front shapes by computing the zero level sets.

As a consequence of discretization, the cost $J(\v p)$ in \eqref{eq:trackingCostH} does not vary in a continuous way, and therefore it is not possible to use gradient-based optimization methods requiring the computation of derivatives of the cost. Without loss of generality, among the various alternatives, in this paper we focus on the mesh adaptive GPS algorithm owing to its computational efficiency and accuracy of the solution \cite{Bertsekas2016NonlinearProgramming}. The optimization routine implementing this method is available in Matlab through the \textit{patternsearch} function of the Global Optimization Toolbox. The tolerance for the stopping criteria is set to $10^{-4}$, while the maximum number of iterations is equal to $2000$. The time integral in the cost function \eqref{eq:trackingCostH} is
discretized over time with a sampling interval $\Delta t$ that varies on a case-by-case basis, as detailed in the following sections. All the tests are carried out on a computer equipped with a 2.5 GHz Intel Xeon CPU with 32 GB of RAM.

\begin{figure*}[tb]
\centering
\subcaptionbox{}{\includegraphics[width=4.4cm,height=4cm]{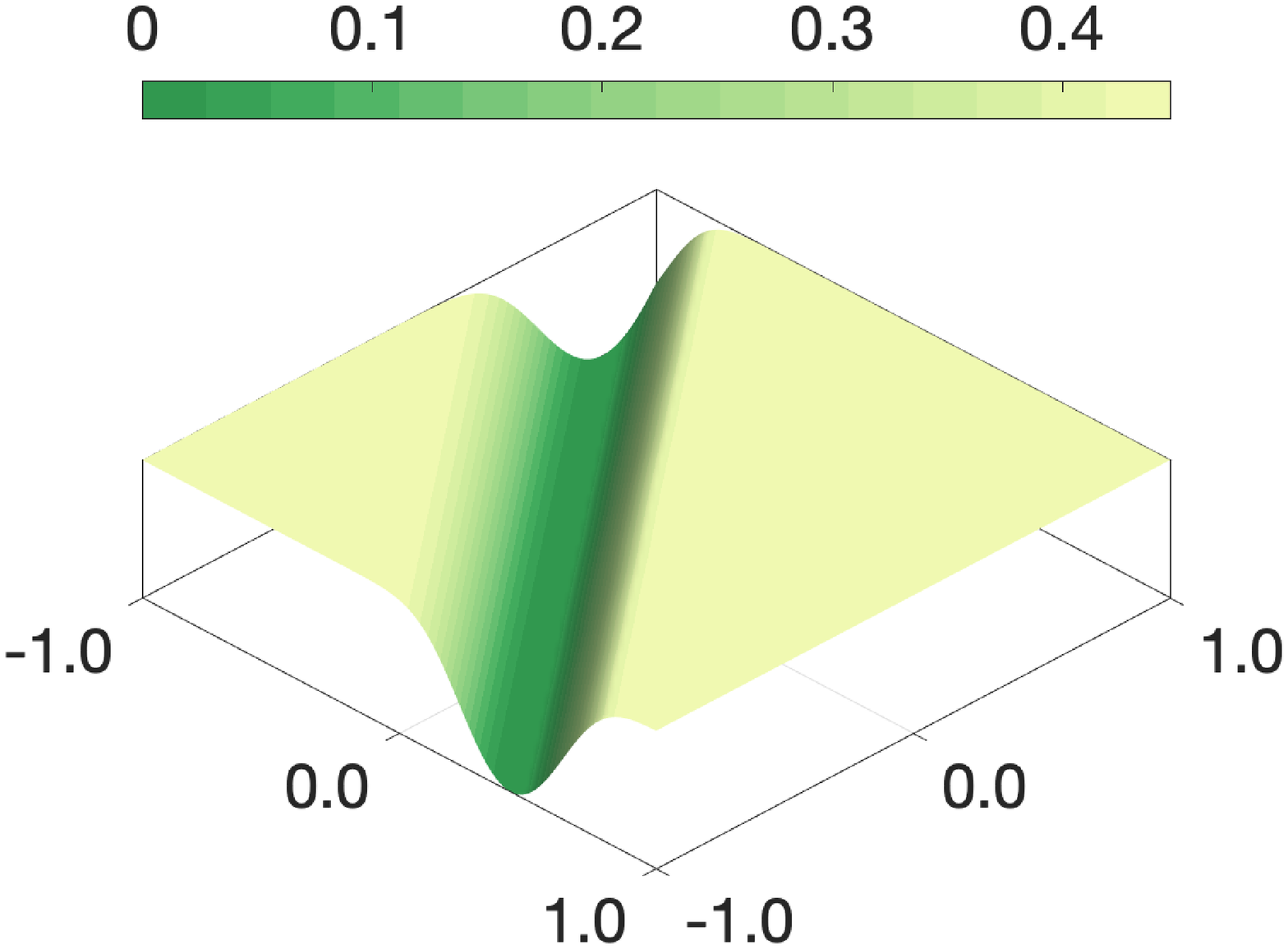}}\label{Valley_Geo.ps}%
%\hfill
\quad
\subcaptionbox{}{\includegraphics[width=4.4cm]{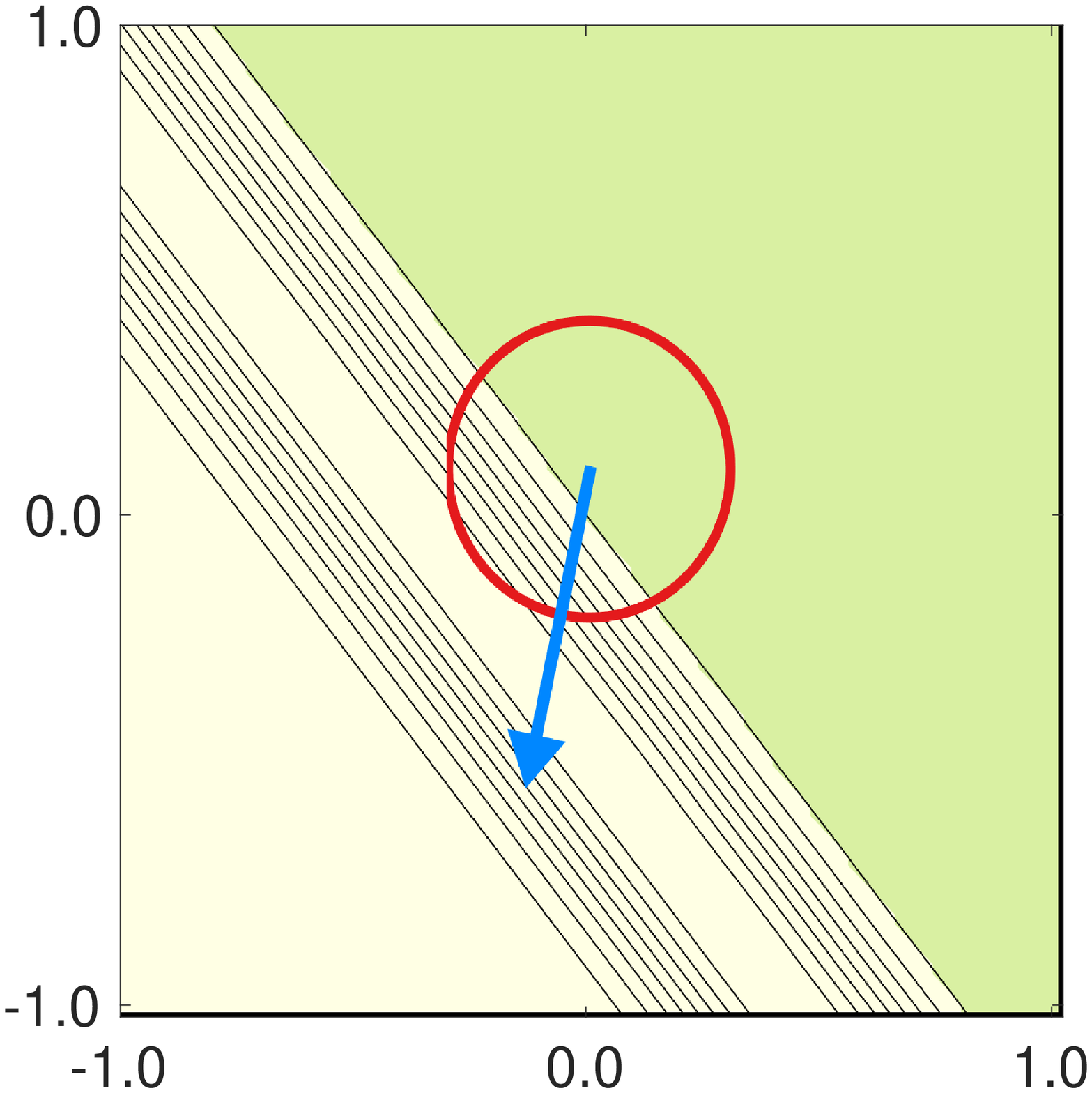}}\label{Valley_Initial.ps}%
\caption{Valley case study: 3D elevation map (a) and initial condition (b). The fire front at $t_0$ is in red, the wind direction is in blue, the elevation level curves are in black, fuel A is in light-yellow, and fuel B is in light-green.}
\label{fig:valley}
\end{figure*}

\begin{figure*}[tb]
\centering
\subcaptionbox{}{\includegraphics[width=4.4cm]{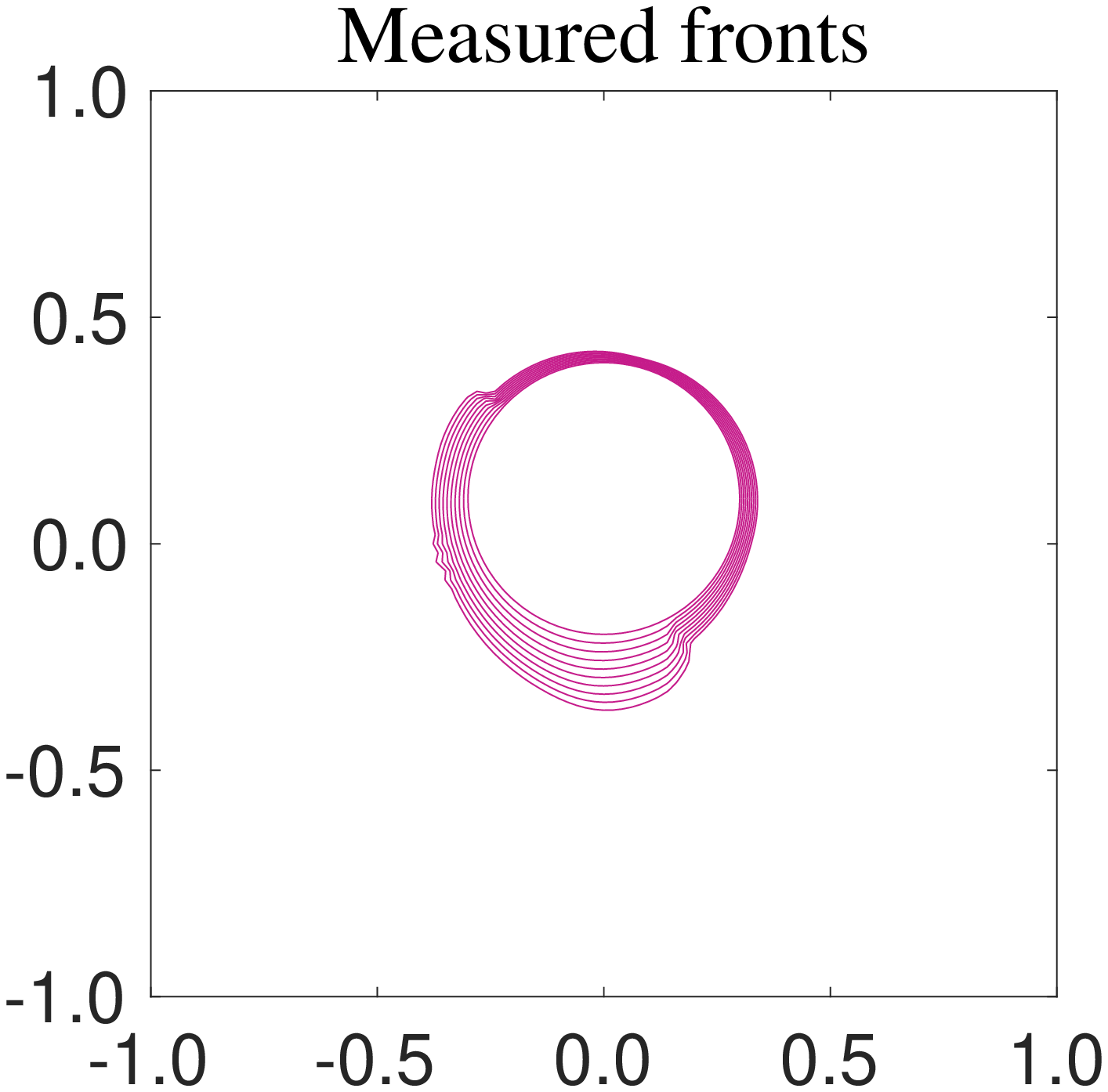}}\label{Valley_Real.ps}%
\subcaptionbox{}{\includegraphics[width=4.4cm]{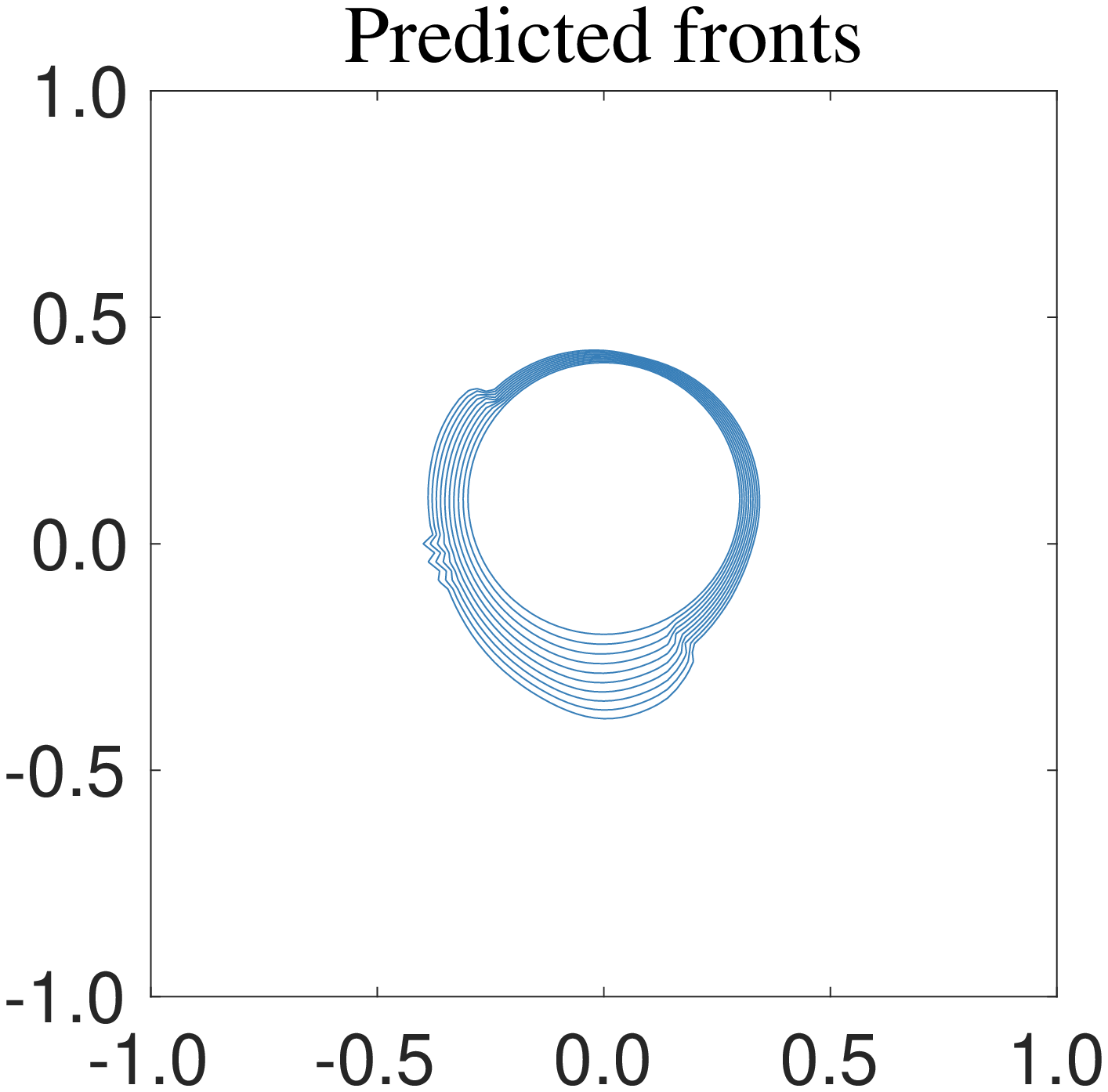}}\label{Valley_White.ps}%
\caption{Valley case study: measured fire fronts obtained with $\v p=\v p^*$ (a) and predicted fire fronts computed with $\v p=\vh p$ (b).}
\label{fig:valleyRealSim}
\end{figure*}

\begin{figure*}[tb]
\centering
\subcaptionbox{}{\includegraphics[width=4.5cm]{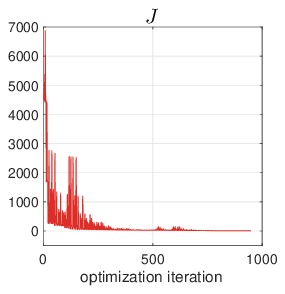}}\label{Valley_Cost.ps}%
\hfill
%\hspace{0.5cm}
\subcaptionbox{}{\includegraphics[width=12.0cm]{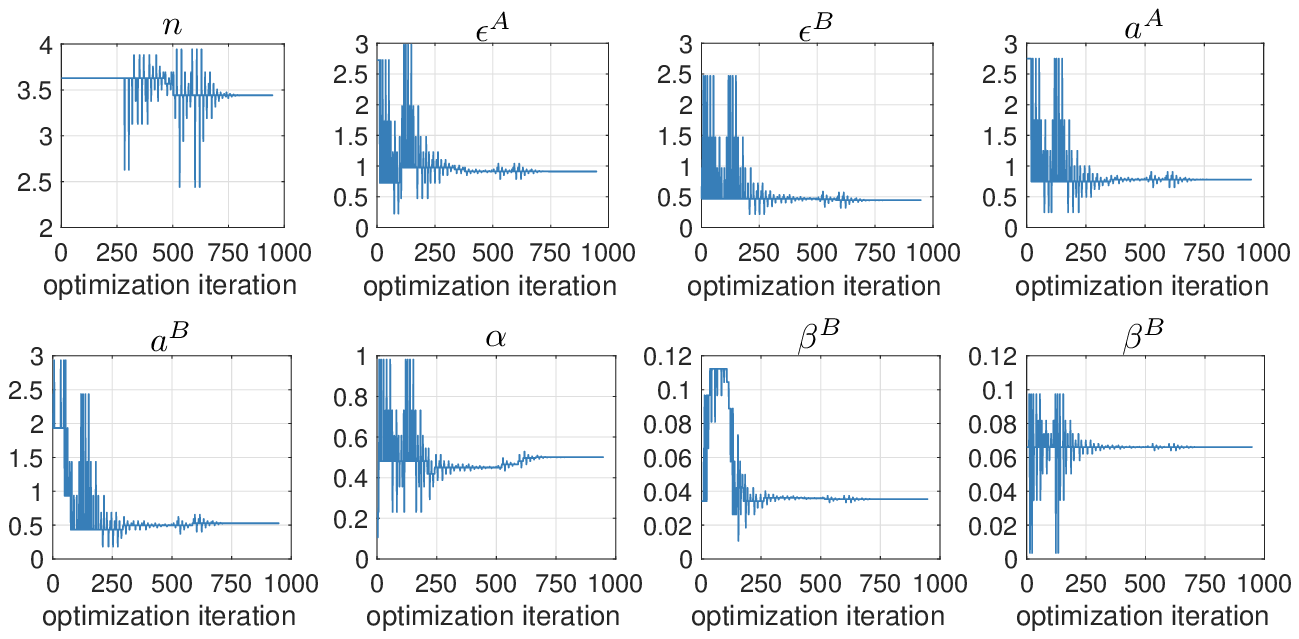}}\label{Valley_Parameters.ps}%
\caption{Valley case study: trends of the cost function (a) and of the parameters (b) during the estimation procedure.}
\label{fig:valleyTrends}
\end{figure*}

\subsection{Model Testing}
\label{sec:sim:testing}

Along the lines of \cite{Lautenberger13} and others, the model described in Section \ref{sec:model} is first tested in four simple cases taken from the reference literature. In all cases, a circle is considered as starting shape of the fire front, and a squared $[-1,1] \times [-1,1]$ domain $\Omega$ is chosen and discretized by means of a $101\times101$ regular grid; the same fuel is considered in all the domain. The  following parameters are adopted: $n=3$, $\epsilon=0.4$, $a=0.5$, $\alpha=0.5$, and $\beta=0.02$. 

In the first case (see Fig. \ref{fig:figsim}(a)), no wind is considered and a final time $t_f$ equal to $0.5$ is fixed starting from $t_0=0$. Owing to the absence of wind, the fire front maintains its circular shape while propagating, as expected \cite{Lo2012AMethod}. In the second case, depicted in Fig. \ref{fig:figsim}(b), a wind with speed components $(-0.5,-0.5)$ is considered over a time span from $t_0=0$ to $t_f=0.5$. The fire front evolves predominantly in the wind direction with an elliptical shape, as shown also in \cite{Mallet2009ModelingMethods, Lo2012AMethod, HuabsomboonImplicitModel}. 

In the third and fourth cases (reported in Fig. \ref{fig:figsim}(c) and Fig. \ref{fig:figsim}(d), respectively), no wind is considered, the landscape elevation is non-constant, and the initial and final times are set equal to $0$ and $1.5$, respectively. In more detail, a hill is present in the third case and, when the fire goes uphill, its spread rate increases (the zero level set curves of $\phi$ considered at a constant time step are more distant one to the others). When the fire is close to the top of the hill, its velocity decreases and it spreads downhill very slowly (the zero level set curves are close one to the others). In the fourth case, a depression is present and the fire behavior is similar, with a lower spread rate when the front moves downhill the depression and a higher one when it moves uphill. Again, such behaviors are in accordance with \cite{Lo2012AMethod, Mallet2009ModelingMethods}. 

Hence, we conclude that the model described in Section \ref{sec:model} behaves as expected and is able to predict future evolution of fire fronts. Thus, we can proceed with the optimal estimation of its parameter to adapt the predicted fronts to the available measures.

\subsection{Parameter Estimation: Valley Case Study}
\label{sec:sim:valley}

In this case study, we consider a region with a terrain characterized by a depression, as shown in Fig. \ref{fig:valley}(a). The initial shape of the fire front is a circle (see Fig. \ref{fig:valley}(b)), and the wind speed components are equal to $(-0.5, -2)$. The domain is given by $\Omega=[-1,1]\times[-1,1]$, and it is discretized by means of a grid made up of $101\times 101$ cells. Two kinds of fuels are present: the fuel A, represented by the light-yellow region in Fig. \ref{fig:valley}(b), is characterized by $\epsilon^A=0.8$, $a^A=0.7$, and $\beta^A=0.03$, while the fuel B, given by the light-green region, is characterized by $\epsilon^B=0.4$, $a^B=0.4$, and $\beta^B=0.08$. The parameters $n$ and $\alpha$ are equal to 3 and 0.5, respectively, on the whole domain. Hence, we define $\v p^*:=(3, 0.8, 0.4, 0.7, 0.4, 0.5, 0.03, 0.08)\in\Re^8$ as the vector of true  parameters. Such true parameters are first used to run a simulation with sampling, initial, and final times equal to $0.01$, $0$, and
$0.1$, respectively. As shown in Fig. \ref{fig:valleyRealSim}(a), the wildland fire propagates faster in the region with fuel A and slows down going downhill of the depression. 

Then, the obtained fire front shapes are saved with a sampling time $\Delta t = 0.01$ and used as measurement fronts $\Gamma^{\rm meas}(t)$ to perform estimation of the parameter vector $\v p$ using the procedure described in Section \ref{sec:ident}. Toward this end, the initial values of the identified parameters vector $\v p$ are randomly chosen between $\v p_{\rm min} = (2,0.1,0.1,0.1,0.1,0.01,0.12,0.12)$ and $\v p_{\rm max} = (4,3,3,3,3,1,0.001,$ $0.001)$. The trends of the components of $\v p$ and of the cost function $J$ during the optimization procedure are shown in Fig. \ref{fig:valleyTrends}. Optimization was performed in $5409$ s. 

It turns out that, after a transient behavior, the cost and the components of the vector $\v p$ converge to stationary values, which represent the optimal estimates. The final estimate of the parameter vector $\v p$ obtained by the optimization algorithm is $\vh p=(3.44, 0.91, 0.44, 0.78, 0.53, 0.5, 0.035, 0.066)$, while the optimal cost $J(\vh p)$ is equal to 13.

\begin{figure*}[tb]
\centering
\subcaptionbox{}{\includegraphics[width=4.4cm,height=4cm]{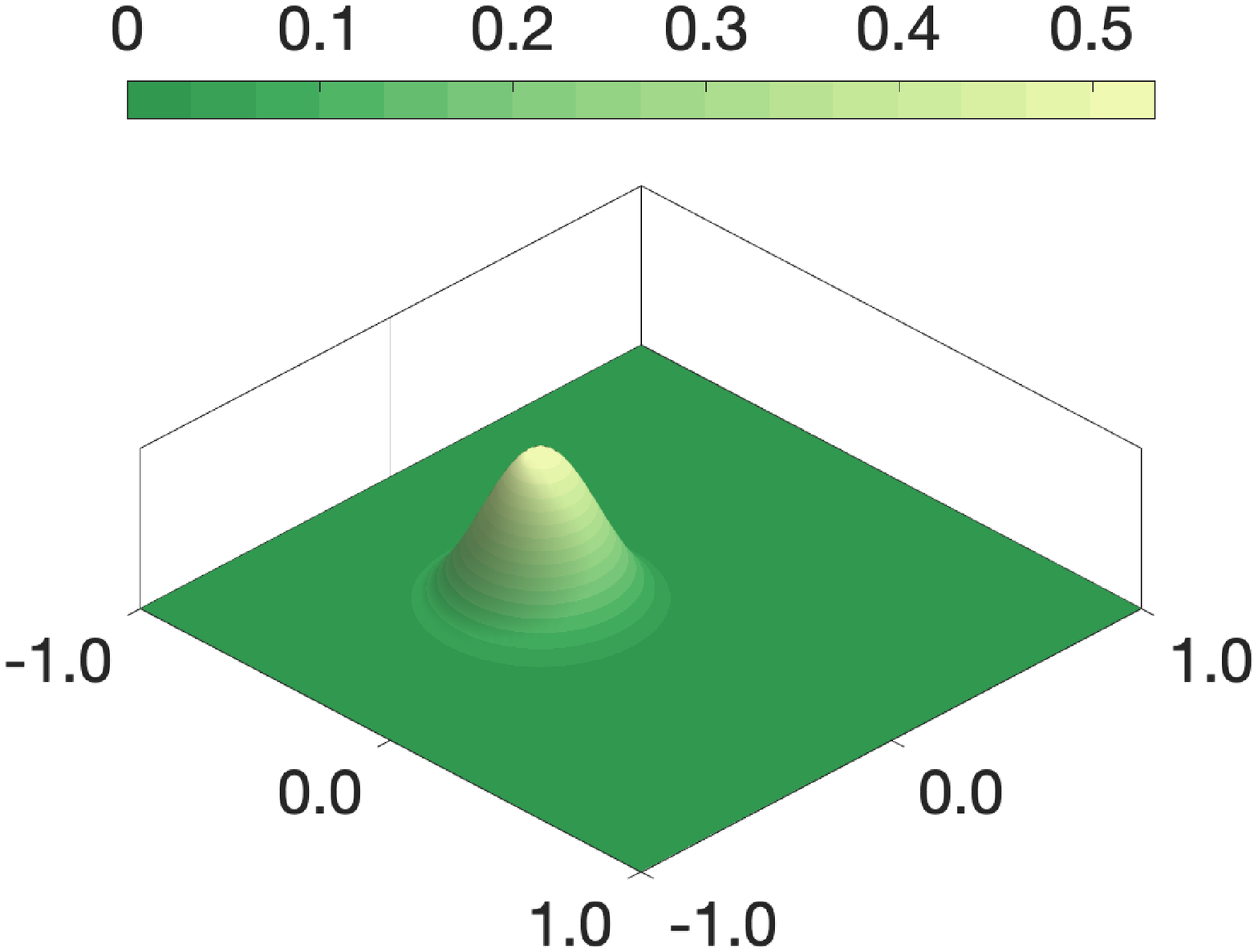}}\label{Hill_Geo.ps}%
%\hfill
\quad
\subcaptionbox{}{\includegraphics[width=4.4cm]{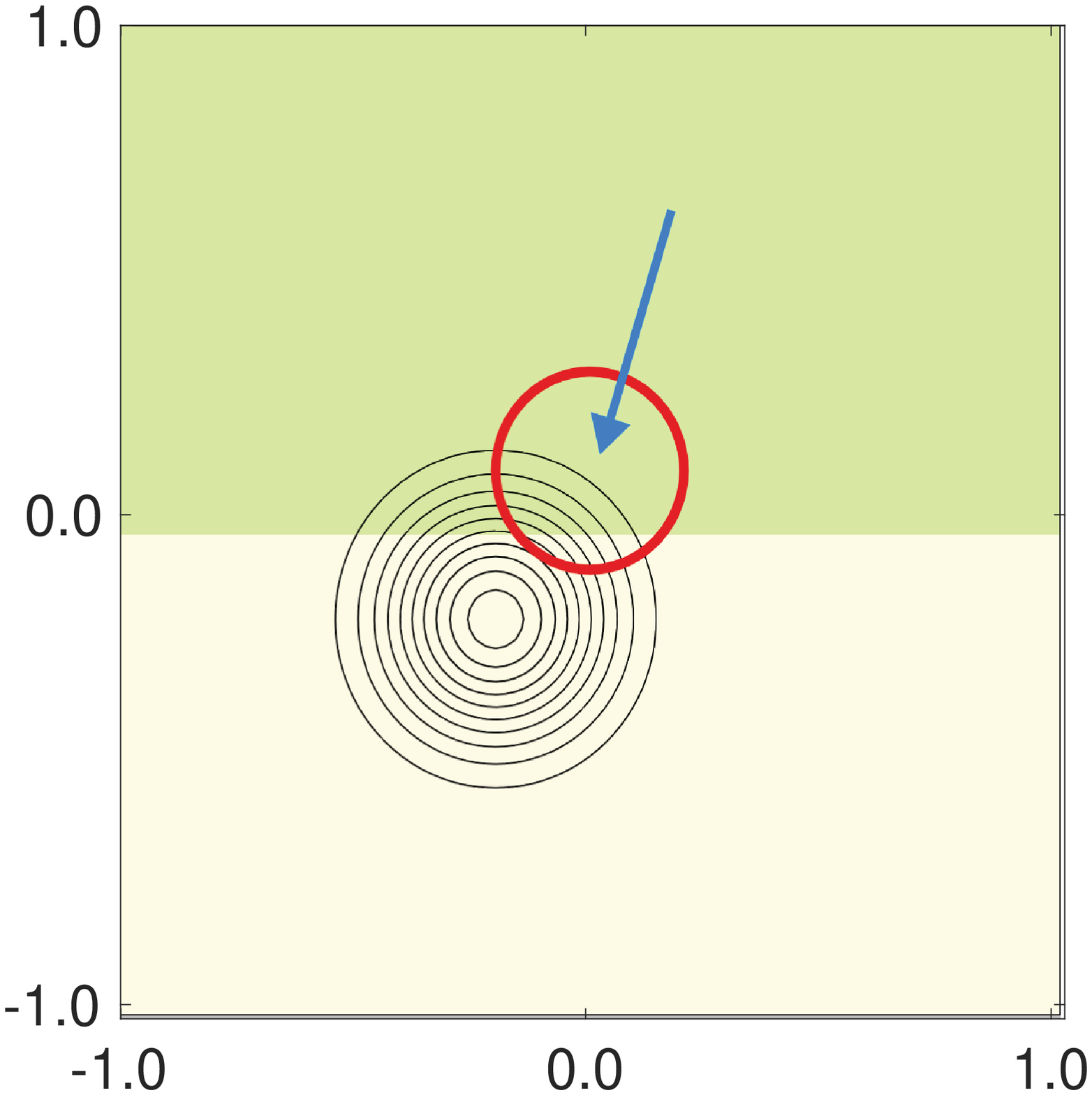}}\label{Hill_Initial.ps}%
\caption{Hill case study: 3D elevation map (a) and initial condition (b). The fire front at $t_0$ is in red, the wind direction is in blue, the elevation level curves are in black, fuel A is in light-yellow, and fuel B is in light-green.}
\label{fig:hill}
\end{figure*}

\begin{figure*}[tb]
\centering
\subcaptionbox{}{\includegraphics[width=4.4cm]{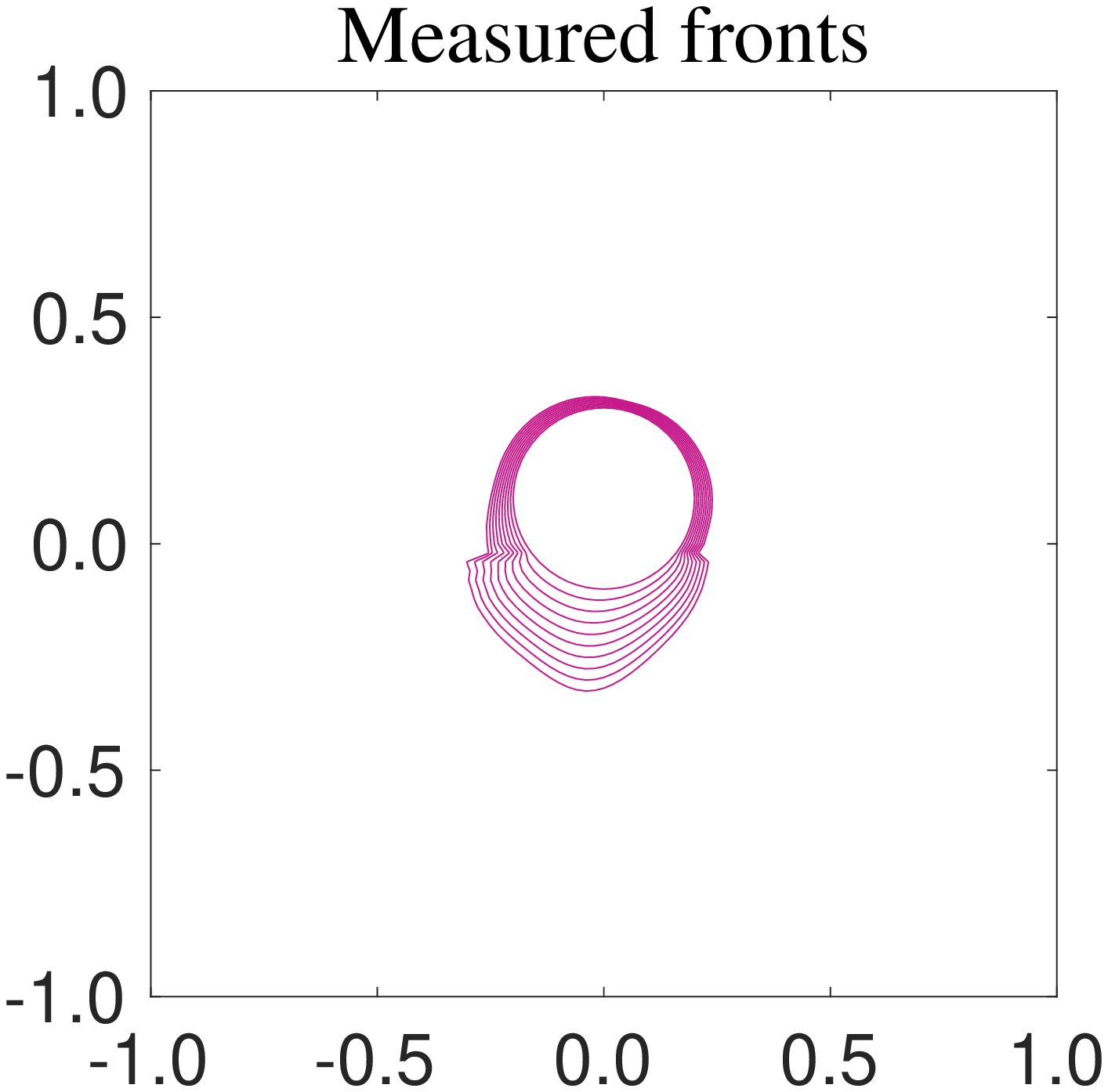}}\label{Hill_Real.ps}%
\subcaptionbox{}{\includegraphics[width=4.4cm]{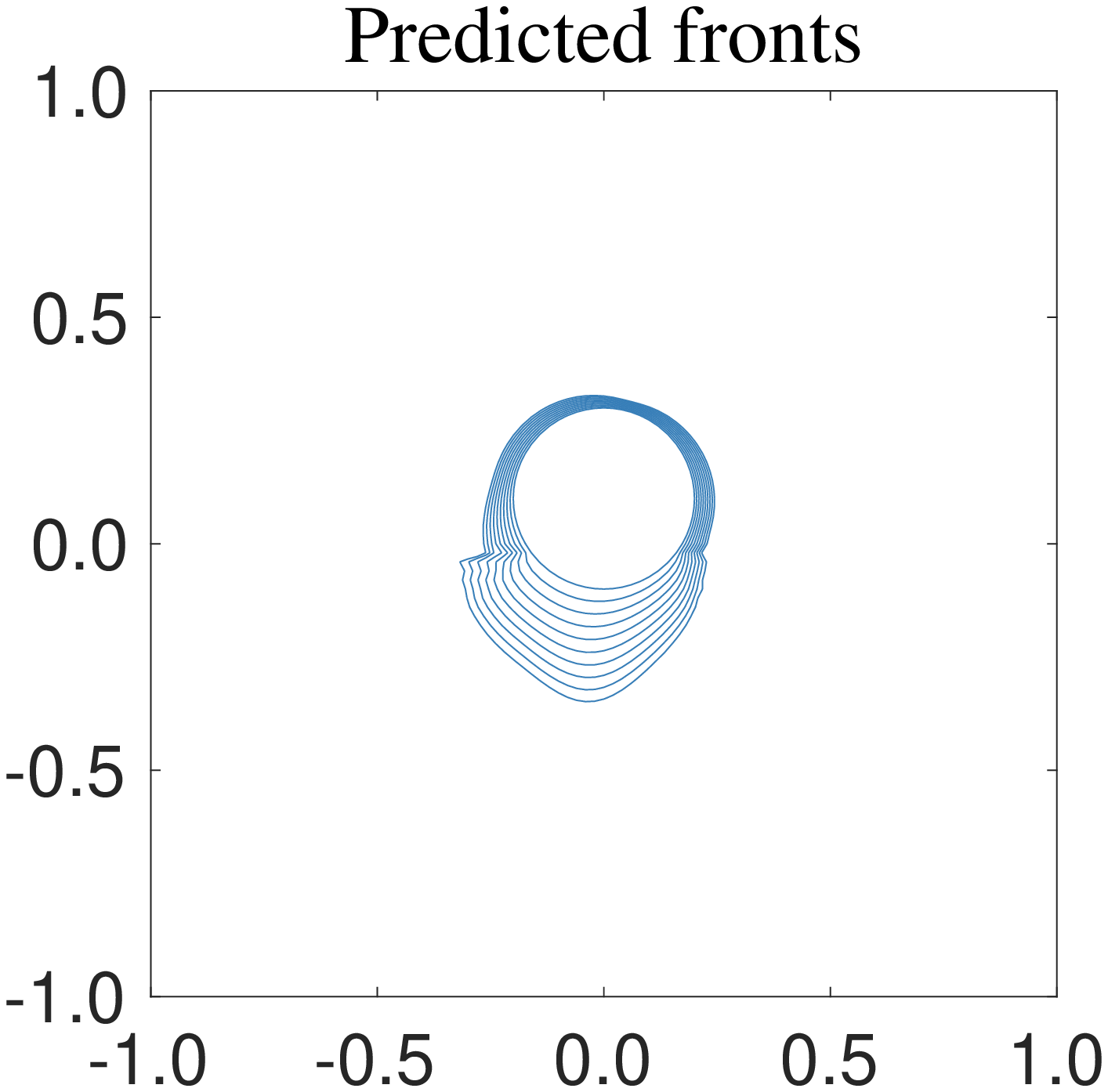}}\label{Hill_White.ps}%
\caption{Hill case study: measured fire fronts obtained with $\v p=\v p^*$ (a) and predicted fire fronts computed with $\v p=\vh p$ (b).}
\label{fig:hillRealSim}
\end{figure*}

\begin{figure*}[tb]
\centering
\subcaptionbox{}{\includegraphics[width=4.5cm]{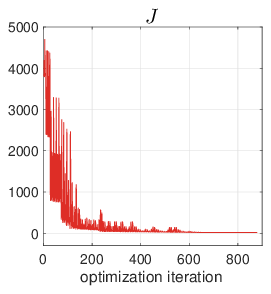}}\label{Hill_Cost.ps}%
\hfill
\subcaptionbox{}{\includegraphics[width=12.0cm]{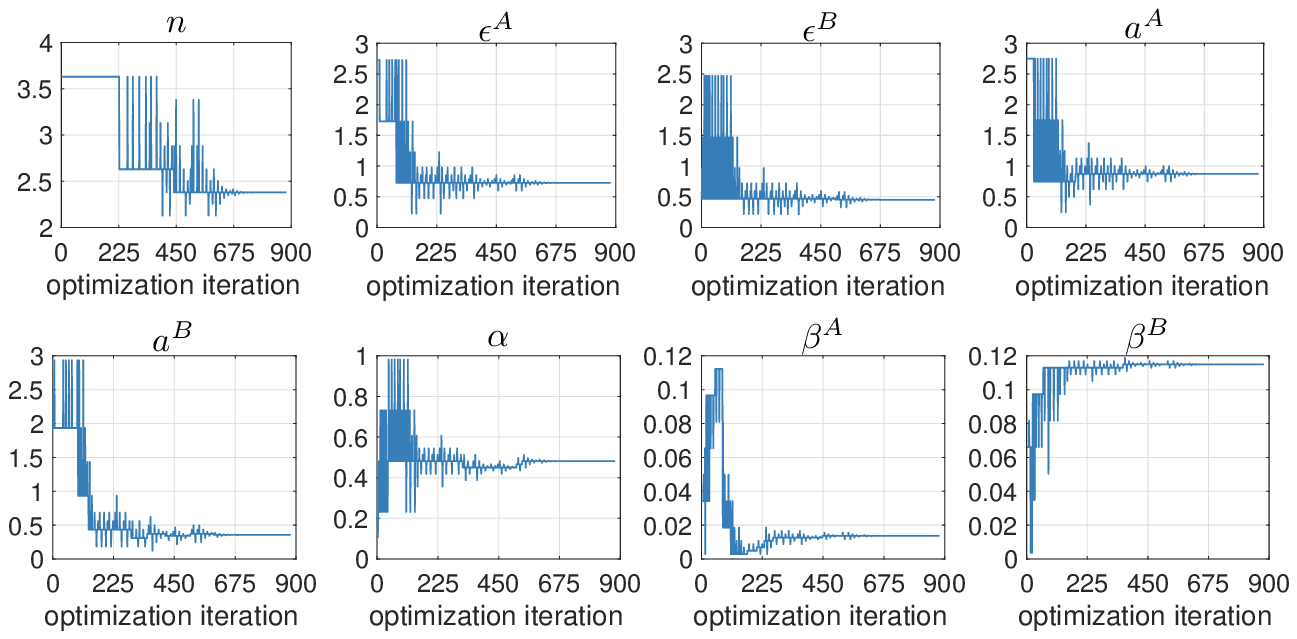}}\label{Hill_Parameters.ps}%
\caption{Hill case study: trends of the cost function (a) and of the parameters (b) during the estimation procedure.}
\label{fig:hillTrends}
\end{figure*}

The accuracy of the estimation can be quantified by computing the following two performance indicators.
\begin{itemize}
	\item[1)] The relative estimation error $e$, defined as the norm of the difference between $\vh p$ and $\v p^*$ with respect to the norm of $\v p^*$, i.e.,
\begin{equation}
e := \frac{\|\vh p - \v p^*\|}{\|\v p^* \|}.
\label{eq:e}
\end{equation}
	\item[2)] The quantity $r$, defined as the ratio between the value of the optimal cost $J(\vh p)$ and the overall number of cells $N$ used to discretize the domain, i.e.,
\begin{equation}
r := \frac{J(\vh p)}{N}.
\label{eq:r}
\end{equation}
It provides the number of cells of the discretization grid where the sign of the functions $\phi^{\rm meas}$ and $\phi$ of the measured and predicted fronts, respectively, differ with respect to the total number of cells.
\end{itemize}
Clearly, the lower are $e$ and $r$, the better is the parameter estimation. 

In the valley case study, we have $e=0.14$ and $r=0.0013$. A simulation of fire propagation with $\v p=\vh p$ provides the fire front snapshot predictions depicted in Fig. \ref{fig:valleyRealSim}(b), which are very close to the corresponding ``measured'' ones reported in Fig. \ref{fig:valleyRealSim}(a). This confirms the effectiveness of the proposed approach for parameter estimation in this case study.

\subsection{Parameter Estimation: Hill Case Study}
\label{sec:sim:hill}

\begin{figure*}[tb]
\subcaptionbox{13:30}{\includegraphics[width=2.15cm]{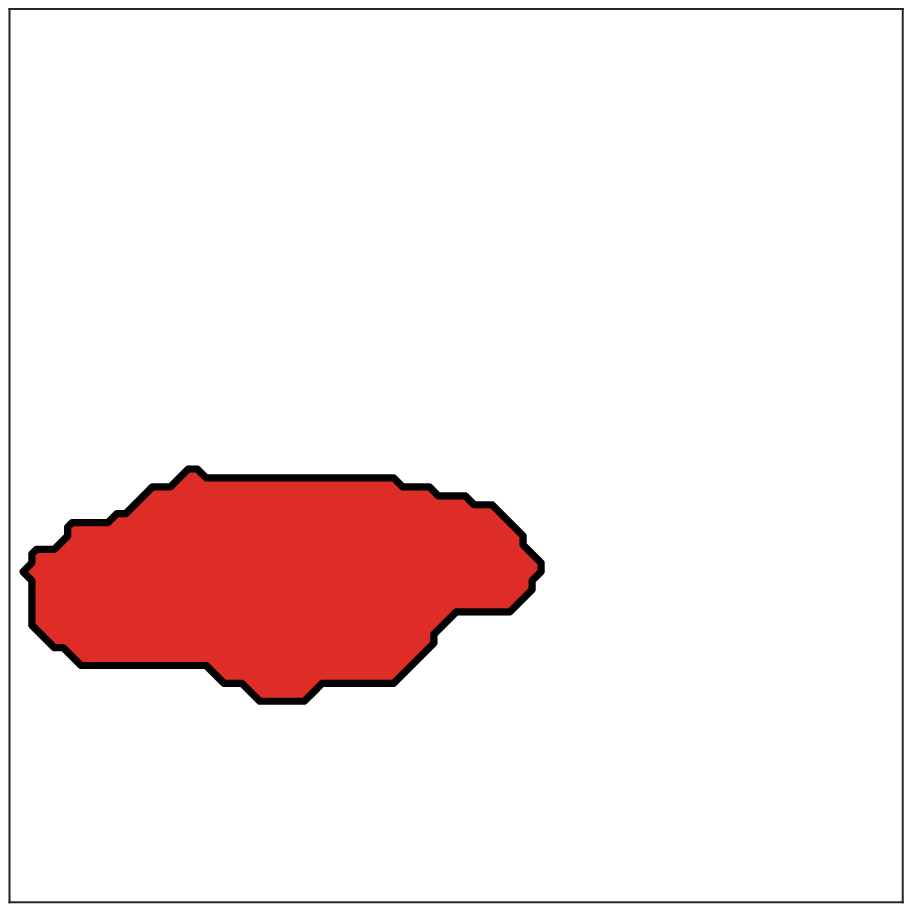}}%
\qquad\quad
\subcaptionbox{13:46}{\includegraphics[width=2.15cm]{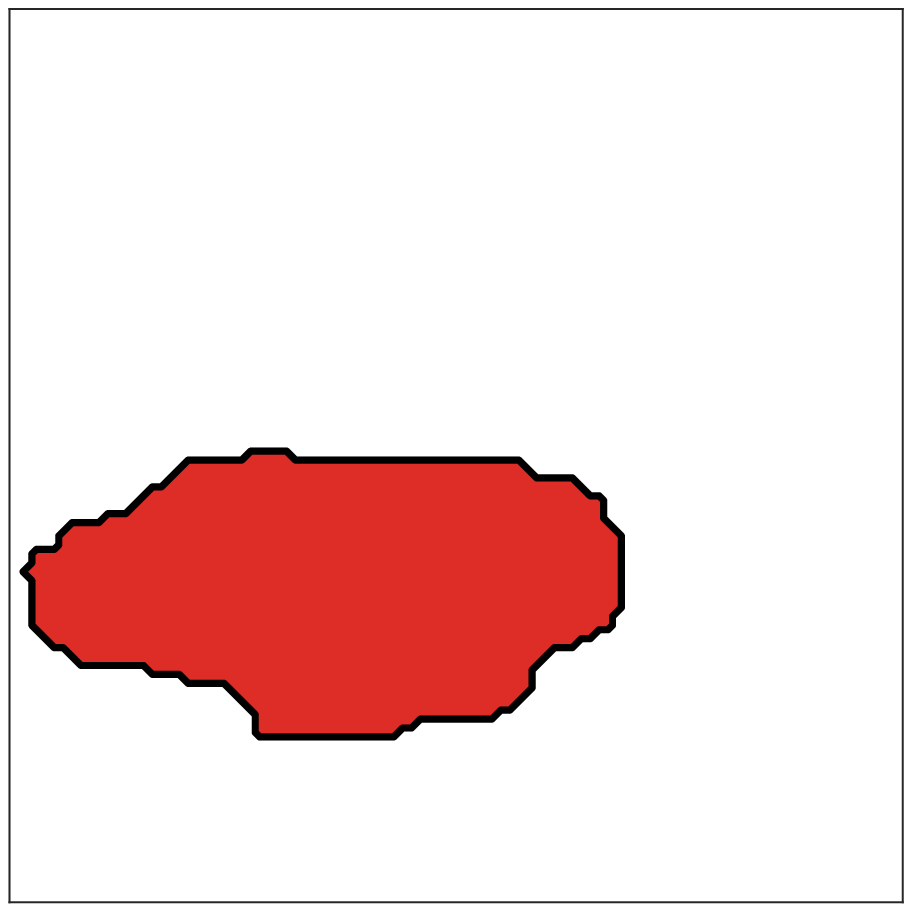}}%
\qquad\quad
\subcaptionbox{14:05}{\includegraphics[width=2.15cm]{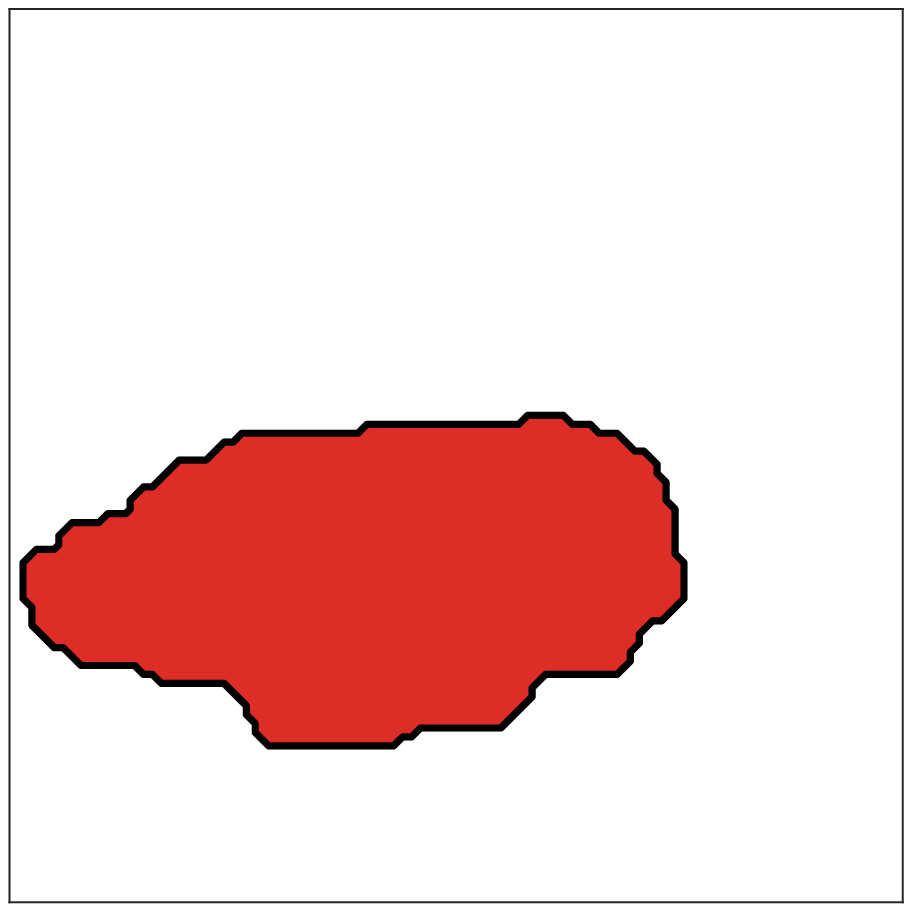}}%
\qquad\quad
\subcaptionbox{14:25}{\includegraphics[width=2.15cm]{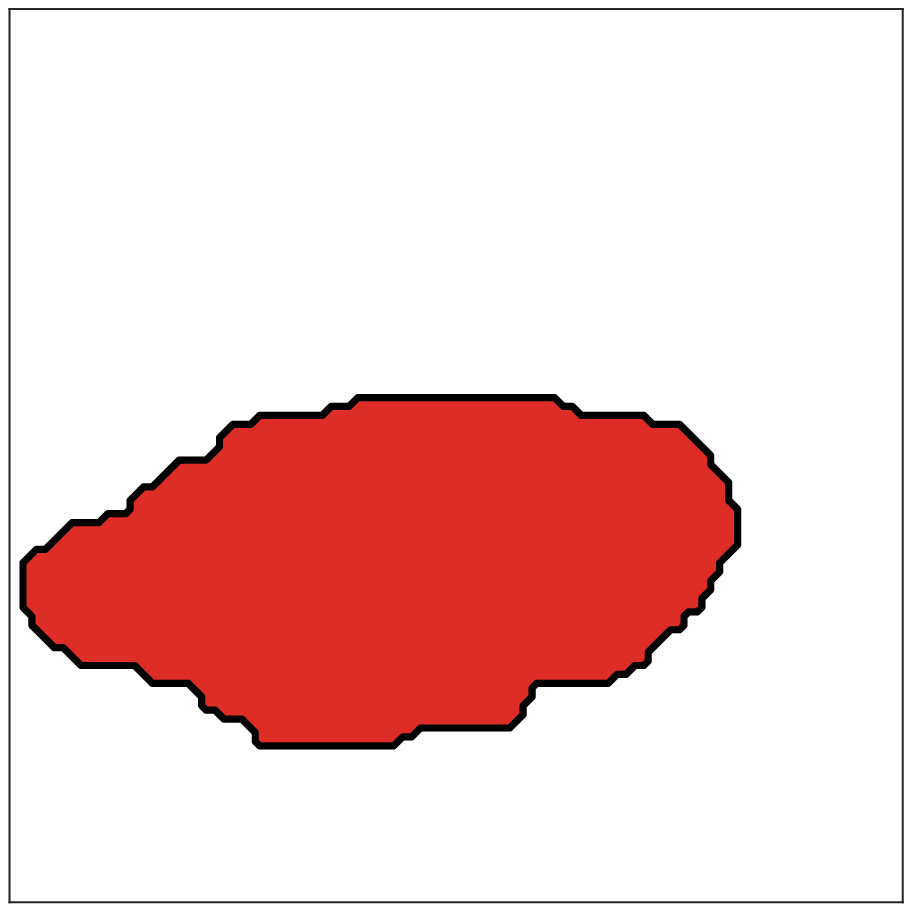}}\\
\\\newline
\subcaptionbox{14:32}{\includegraphics[width=2.15cm]{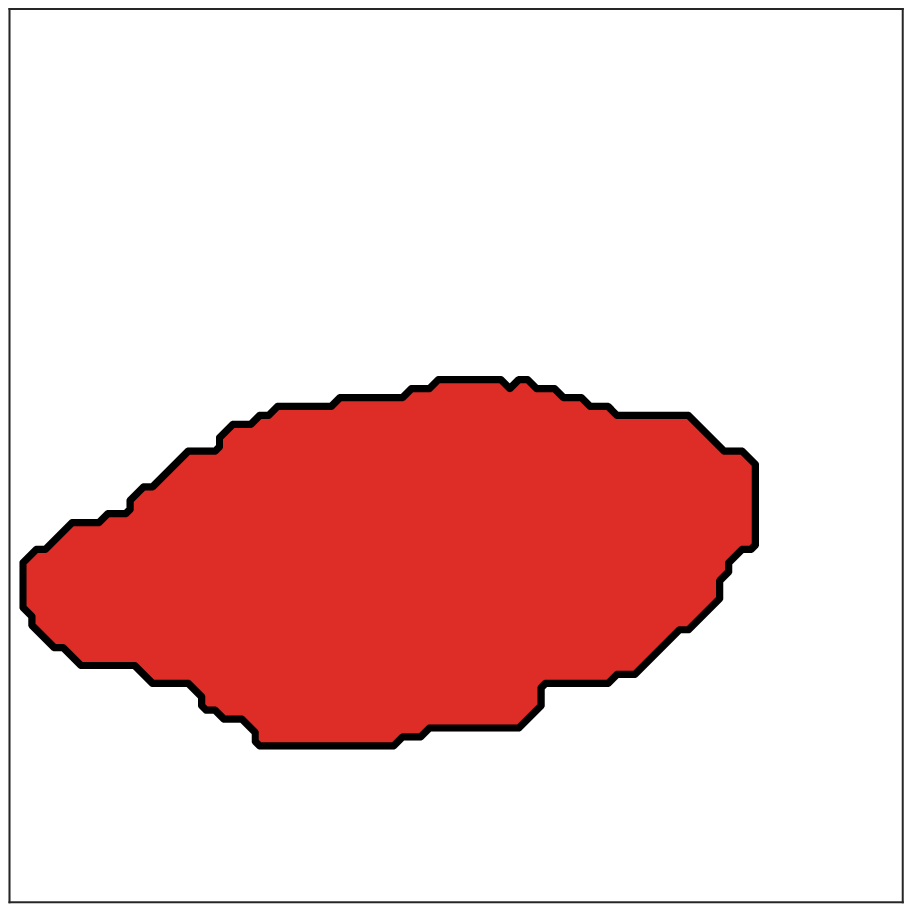}}%
\qquad\quad
\subcaptionbox{14:39}{\includegraphics[width=2.15cm]{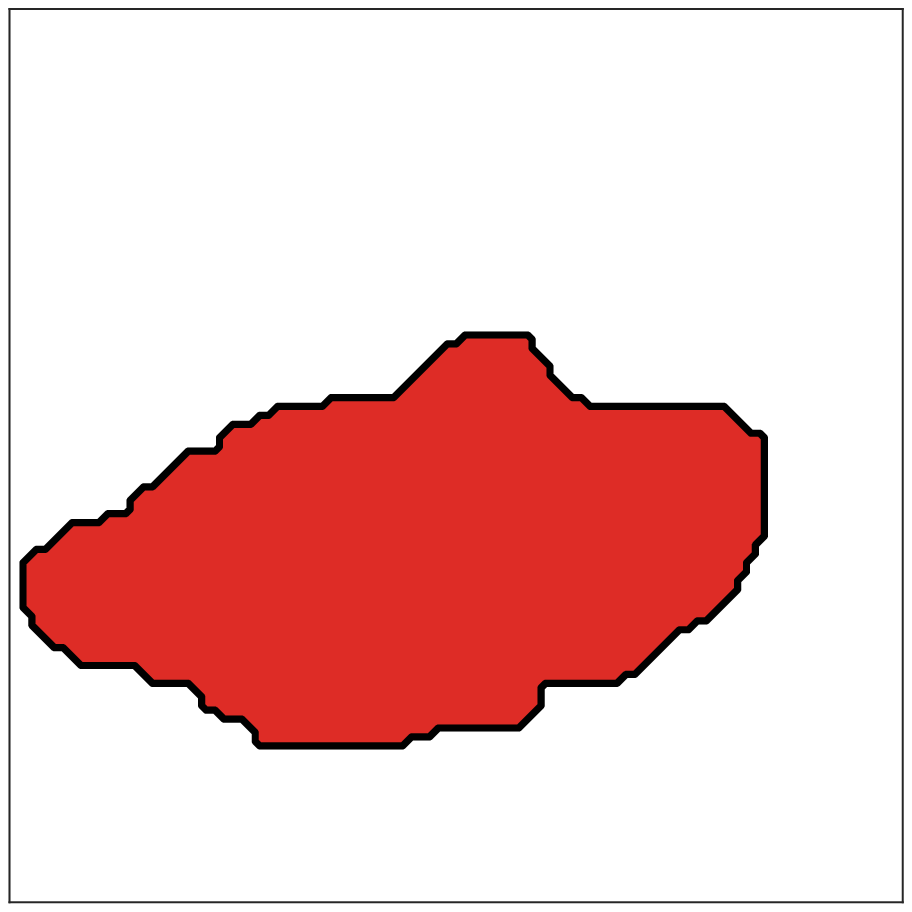}}%
\qquad\quad
\subcaptionbox{14:46}{\includegraphics[width=2.15cm]{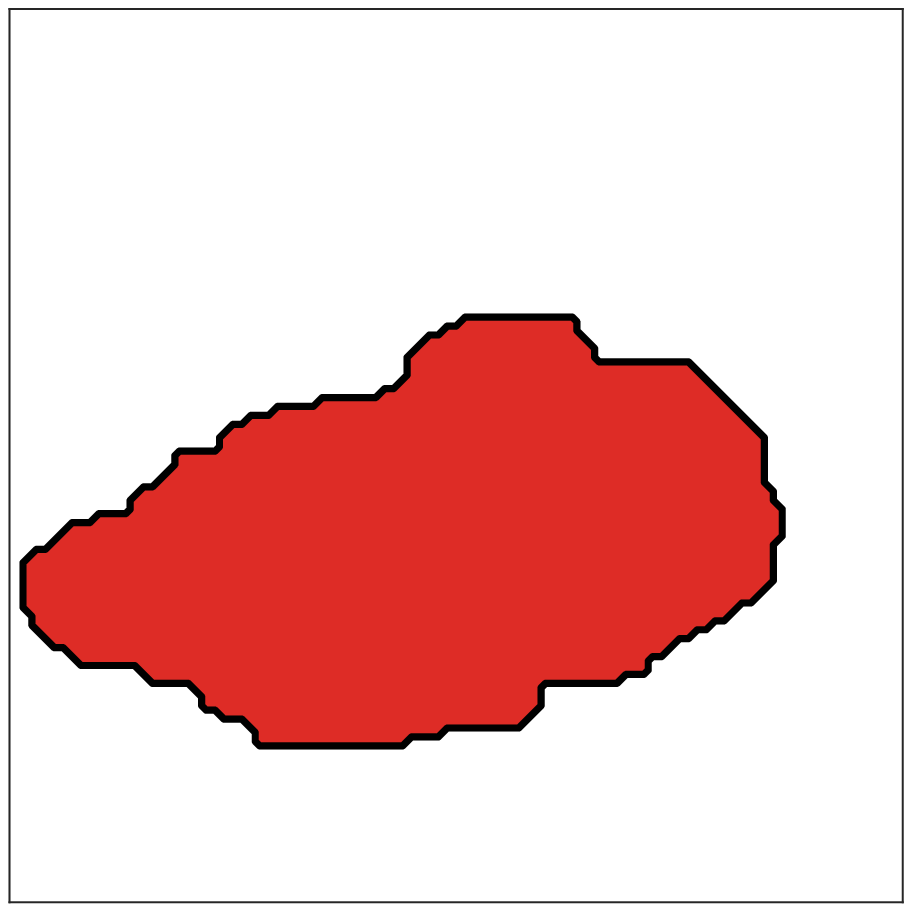}}%
\qquad\quad
\subcaptionbox{14:52}{\includegraphics[width=2.15cm]{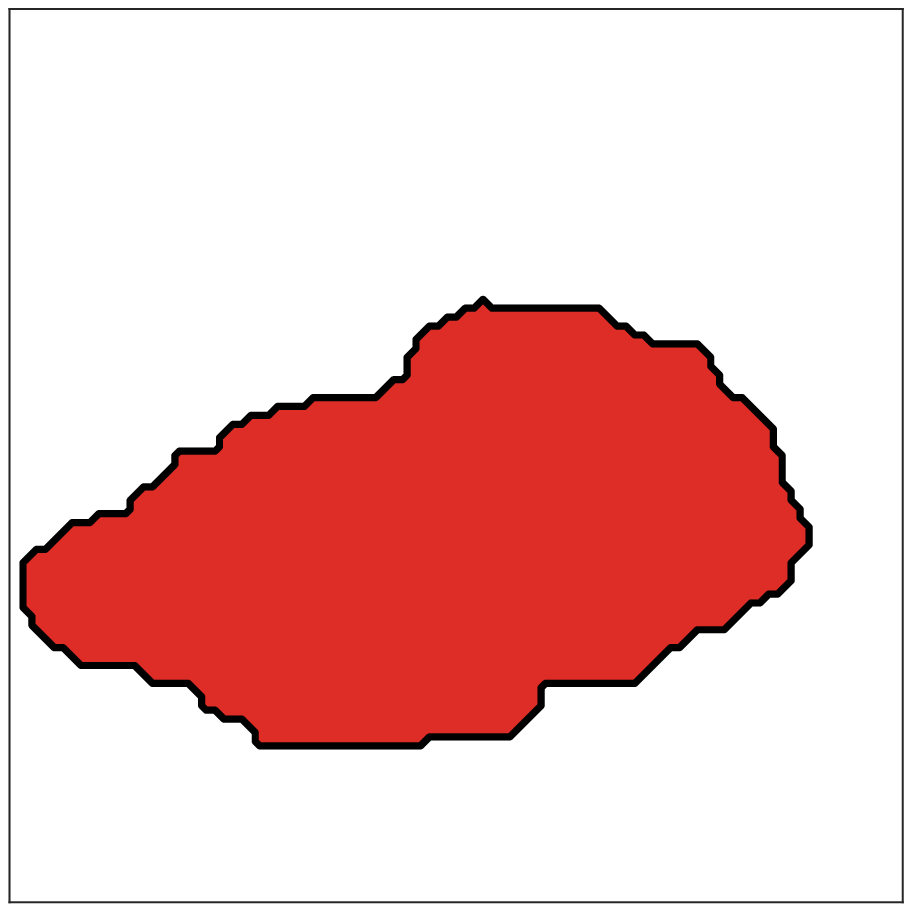}}\\
\\\newline
\subcaptionbox{15:01}{\includegraphics[width=2.15cm]{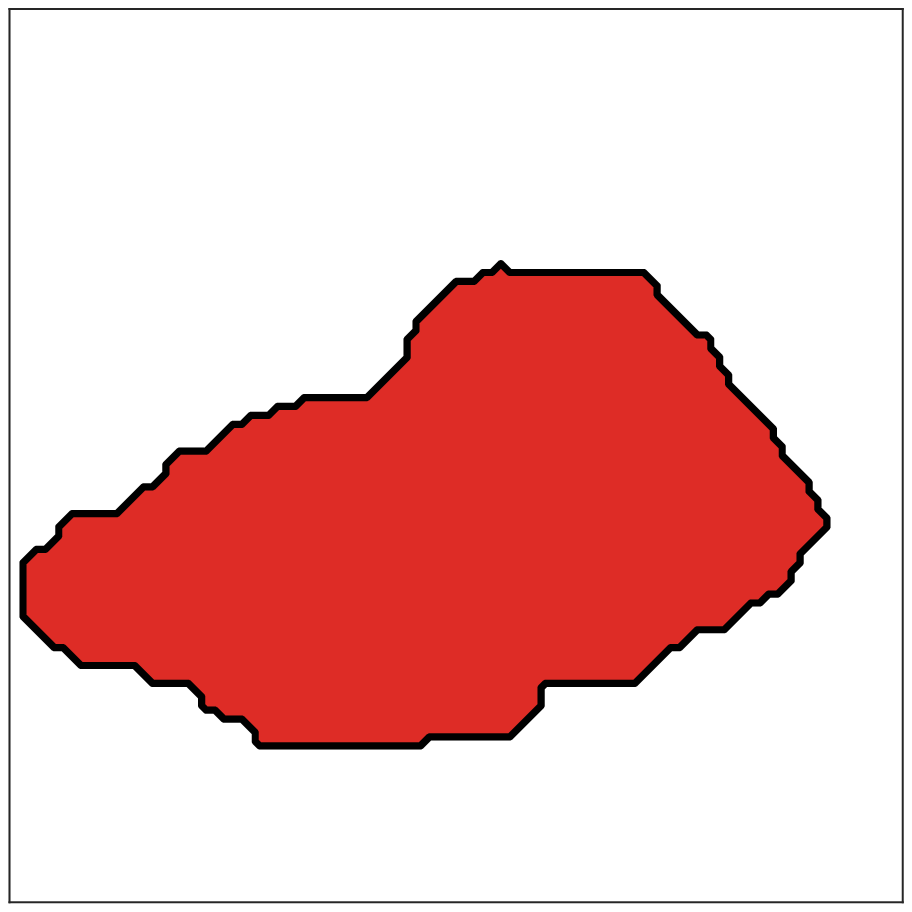}}%
\qquad\quad
\subcaptionbox{15:14}{\includegraphics[width=2.15cm]{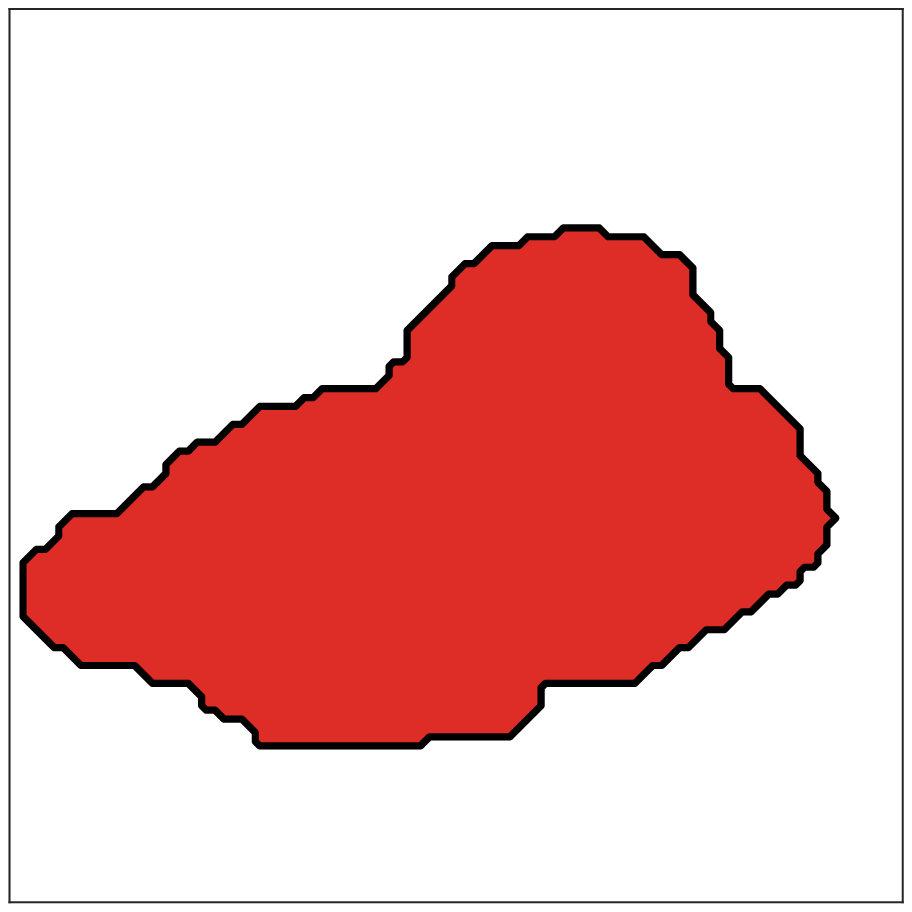}}%
\qquad\quad
\subcaptionbox{15:21}{\includegraphics[width=2.15cm]{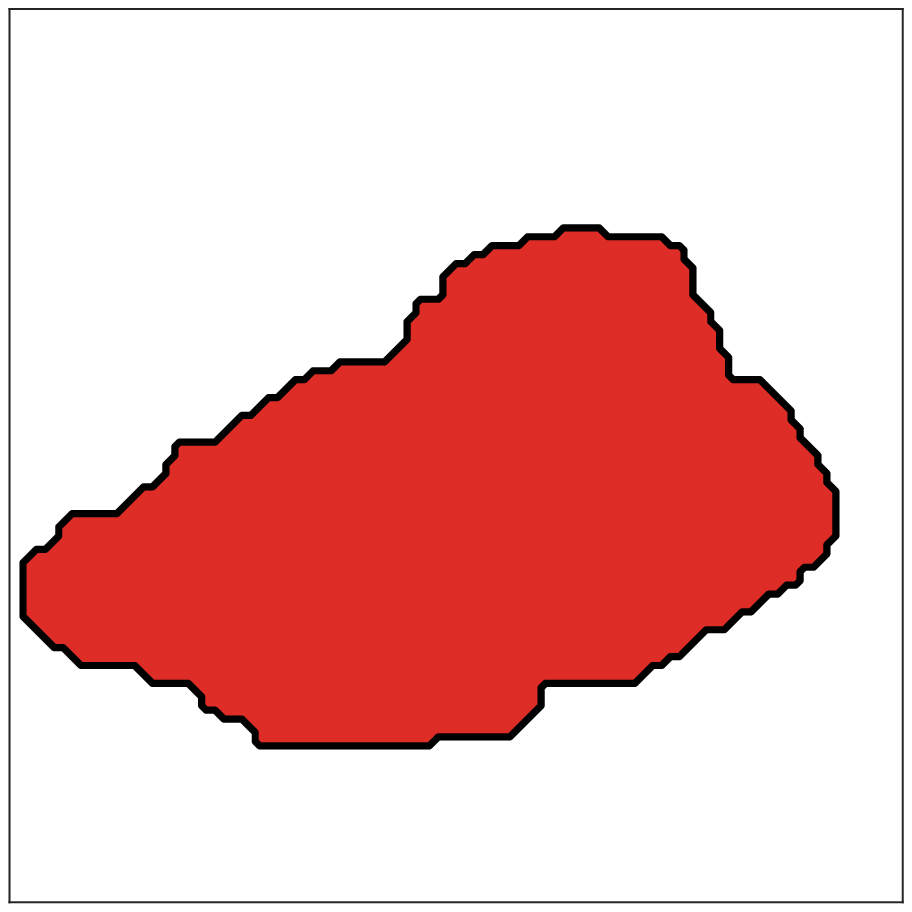}}%
\qquad\quad
\subcaptionbox{15:28}{\includegraphics[width=2.15cm]{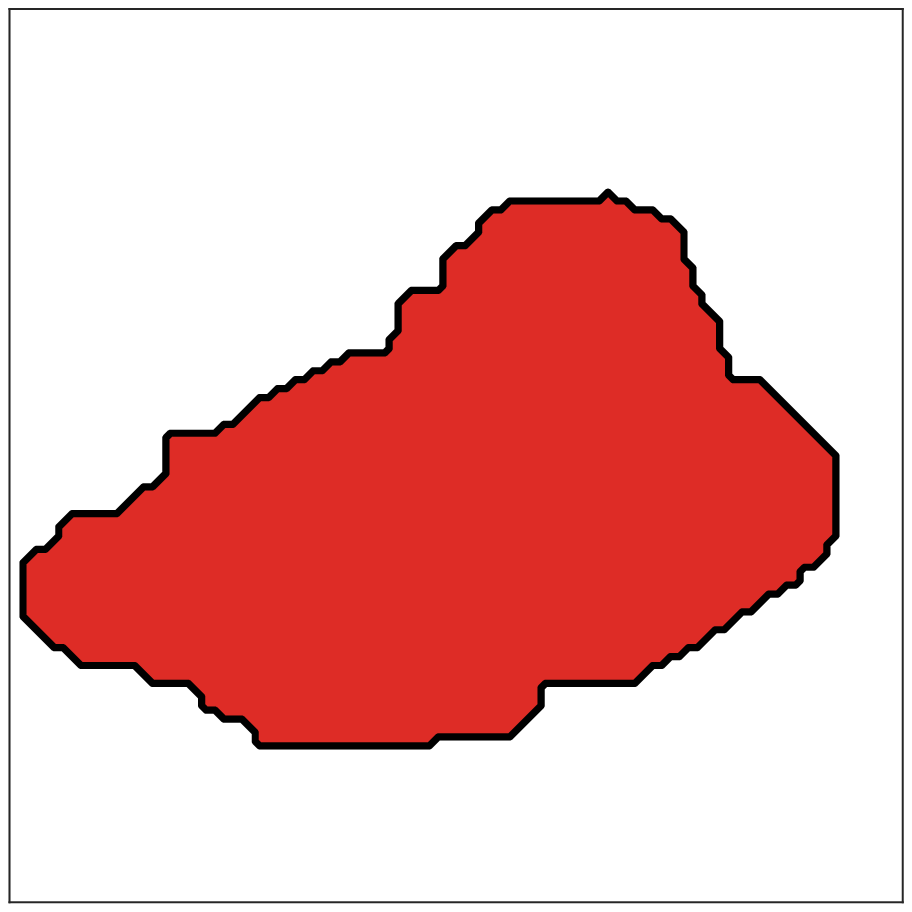}}\\
\\\newline
\subcaptionbox{15:35}{\includegraphics[width=2.15cm]{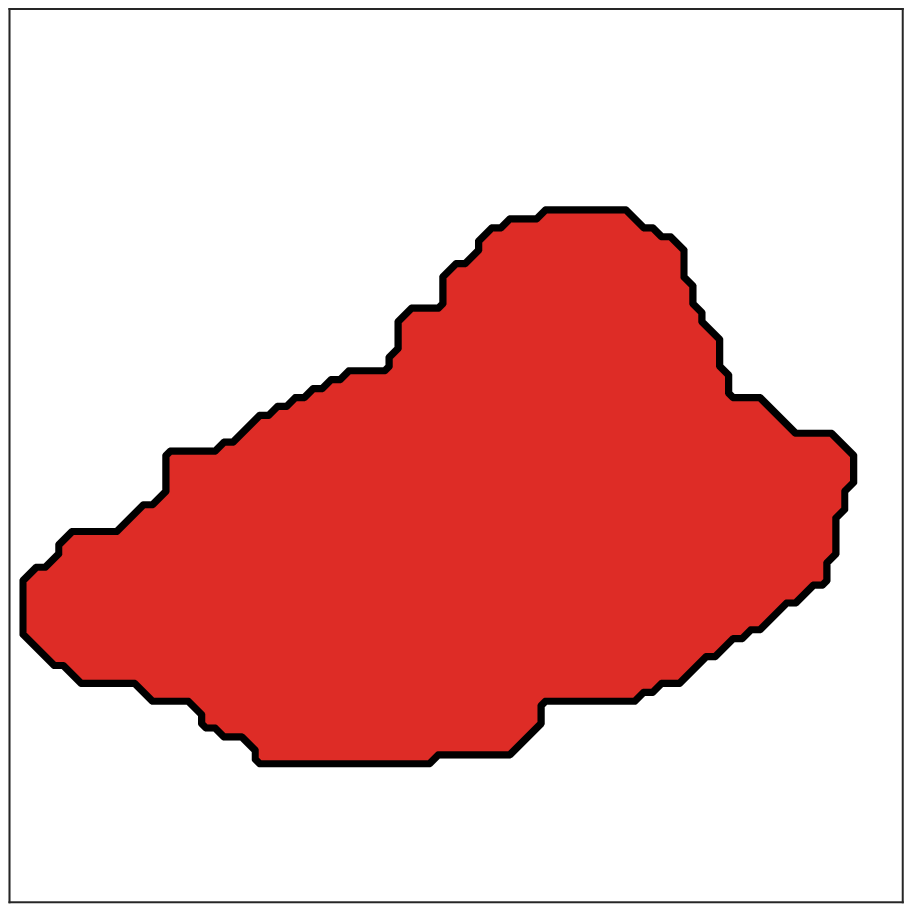}}%
\qquad\quad
\subcaptionbox{15:42}{\includegraphics[width=2.15cm]{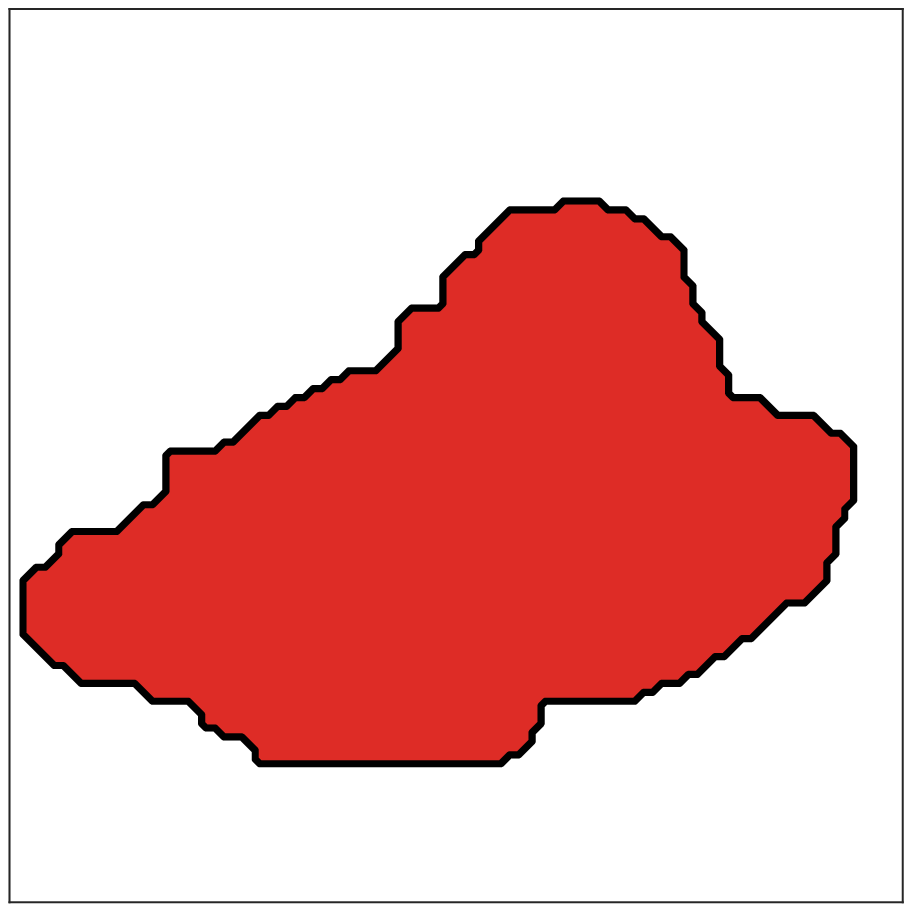}}%
\qquad\quad
\subcaptionbox{15:48}{\includegraphics[width=2.15cm]{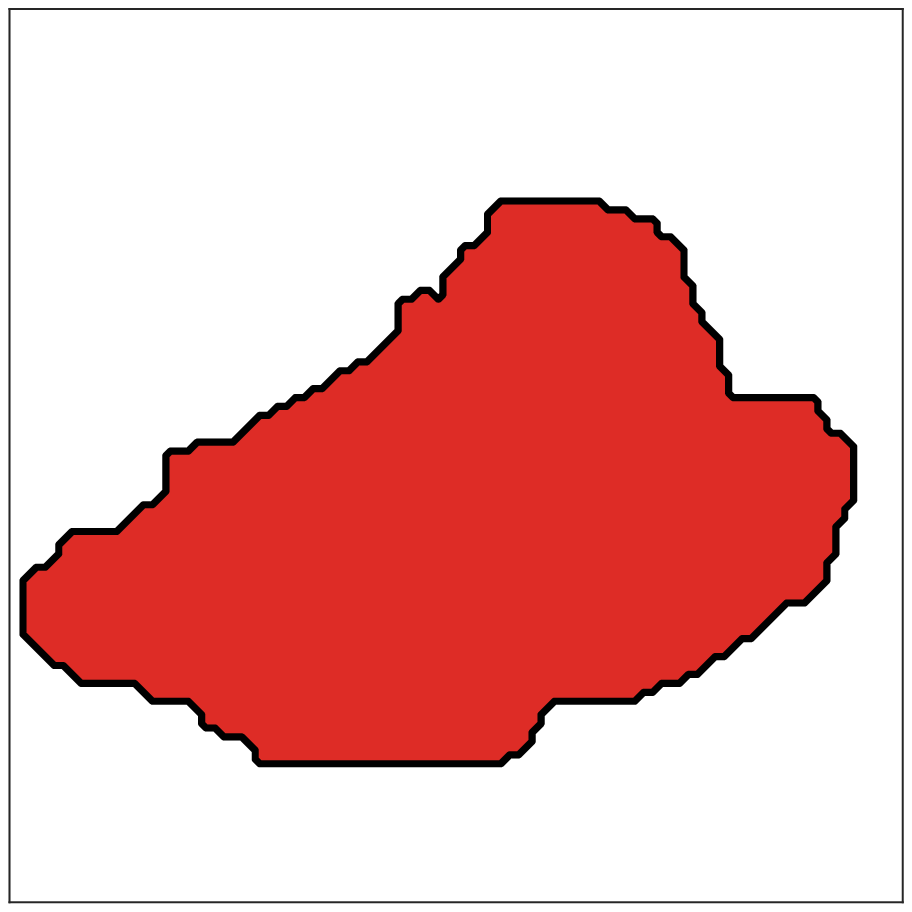}}%
\qquad\quad
\subcaptionbox{15:55}{\includegraphics[width=2.15cm]{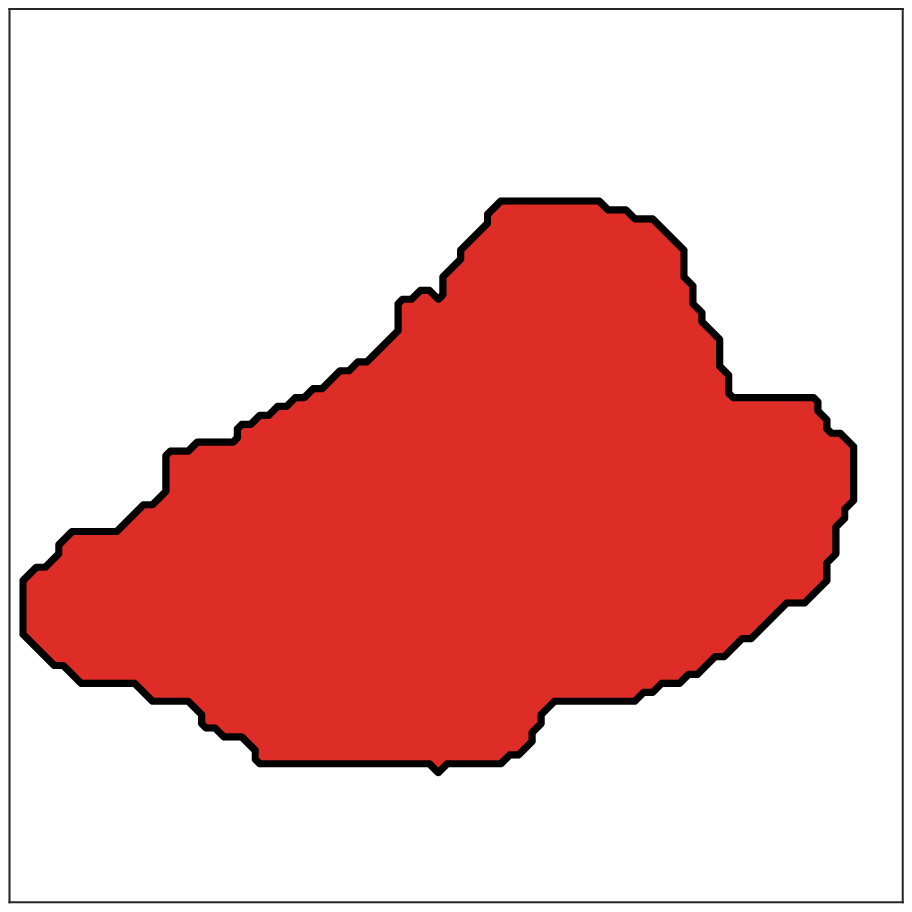}}%
\caption{Evolution of the Troy fire between 13:30 and 15:55. The red area represents the burnt region, while the black contour is the fire front.}
\label{fig:troySatellite}
\end{figure*}

This case study is characterized by the presence of a hill instead of a valley in the domain $\Omega$, as showcased in Fig. \ref{fig:valley}(a). The initial shape of the fire front is still a circle (see Fig. \ref{fig:valley}(b)), and the wind speed components are fixed to $(-1,-3)$. The domain is given again by $\Omega=[-1,1]\times[-1,1]$, and it is discretized with a $101\times 101$ grid. Two fuels are present, referred to as A and B, as in the previous case study. The fuel A is in the light-yellow region in Fig. \ref{fig:hill}(b), and it is characterized by $\epsilon^A=0.8$, $a^A=0.6$, and $\beta^A=0.02$, while the fuel B is in the light-green region, and we have
$\epsilon^B=0.5$, $a^B=0.4$, and $\beta^B=0.06$. On the whole domain, $n$ and $\alpha$ are set equal to $3$ and $0.5$, respectively. Thus, the vector of true parameters is given by $\v p^*:=(3,0.8,0.5,0.6,0.4,0.5,0.02,0.06)\in\Re^8$. As in the previous case study, these parameters are first used to run a simulation with sampling, initial, and final time equal to $0.01$, $0$, and $0.1$, respectively. As shown in Fig. \ref{fig:hillRealSim}(a), the fire propagates faster in the region where fuel A is present, it goes faster going uphill and slows down going downhill.

The obtained fire front shapes are saved with a sample time $\Delta t = 0.01$ and used as measurement fronts $\Gamma^{\rm meas}(t)$ to perform estimation of the parameter vector $\v p$ using the procedure described in Section \ref{sec:ident}. The initial values of the parameter vector $\v p$ are randomly chosen between $\v p_{\rm min} = (2,0.1,0.1,0.1,0.1,0.1,0.12,$ $0.12)$ and $\v p_{\rm max} = (4,3,3,3,3,1,0.001,0.001)$. The trends of the components of $\v p$ and of the cost function $J$ during the optimization procedure are shown in Fig. \ref{fig:hillTrends}. 

Likewise in the valley case study, after a transient behavior, the cost and the components of the vector $\v p$ converge to stationary values, which represent the optimal estimates. The final estimate of the parameter vector $\v p$ obtained by the optimization algorithm is $\vh p=(2.38, 0.73, 0.45, 0.87, 0.36,$ $ 0.48, 0.014, 0.11)$, while the optimal cost $J(\vh p)$ is equal to 17. Optimization was performed in $5165$ s. 

The fire front with $\v p=\vh p$ propagates as showcased in Fig. \ref{fig:hillRealSim}(b). Also for this case study, there is a good correspondence between the simulation with real parameters (reported in Fig. \ref{fig:hillRealSim}(a)) and the one run with the identified parameters (see Fig. \ref{fig:hillRealSim}(b)), thus confirming the effectiveness of the estimation procedure. In fact, we have obtained an estimation error $e=0.20$ and a ratio between the value of the optimal cost and the overall number of cells used to discretize the domain equal to $r=0.017$.

\subsection{Parameter Estimation: Troy Case Study}
\label{sec:sims:troy}

\begin{figure*}[tb]
\centering
\subcaptionbox{}{\includegraphics[width=4.4cm,height=4cm]{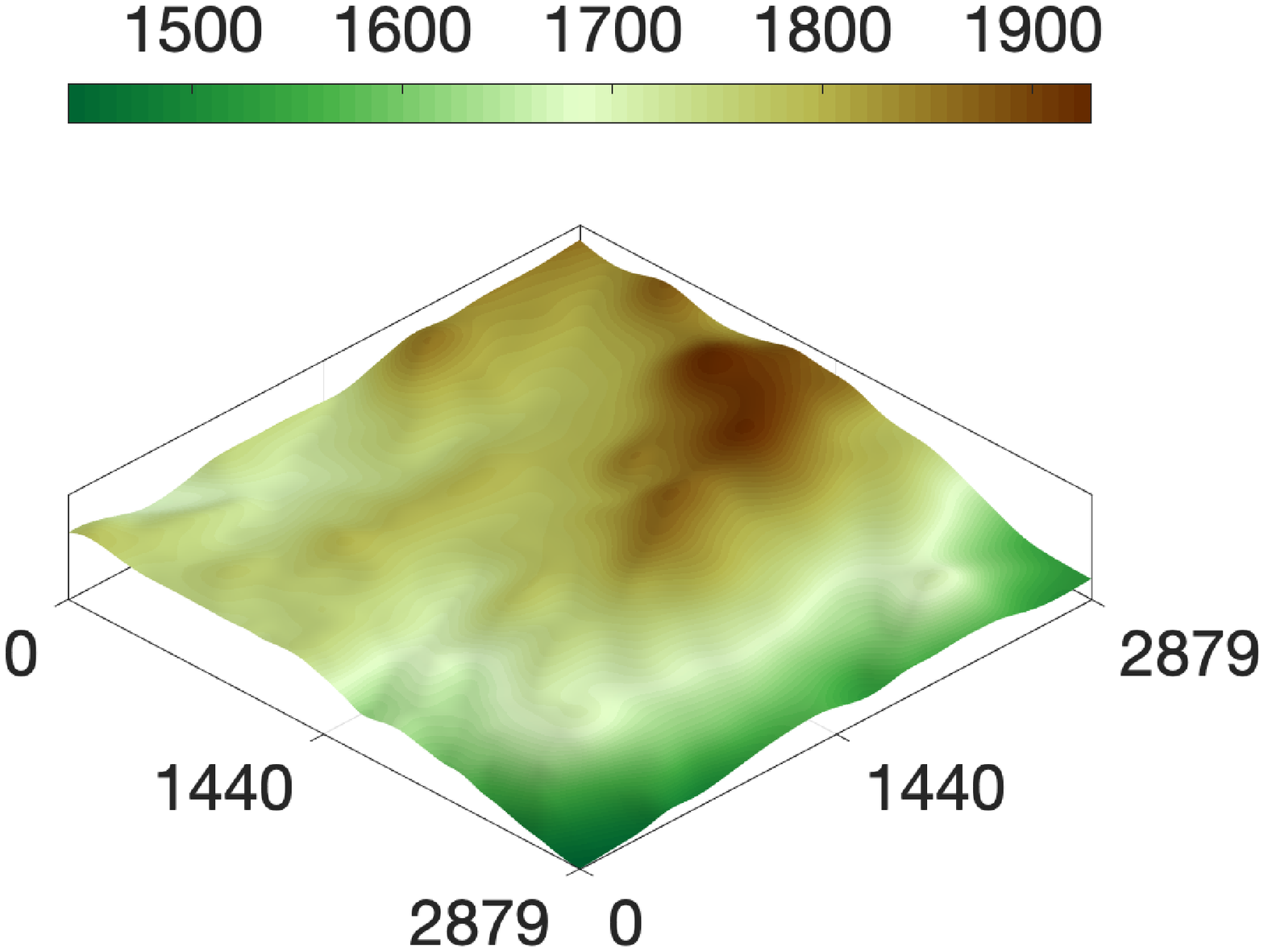}}\label{Troy_Geo.ps}%
%\hfill
\quad
\subcaptionbox{}{\includegraphics[width=4.4cm]{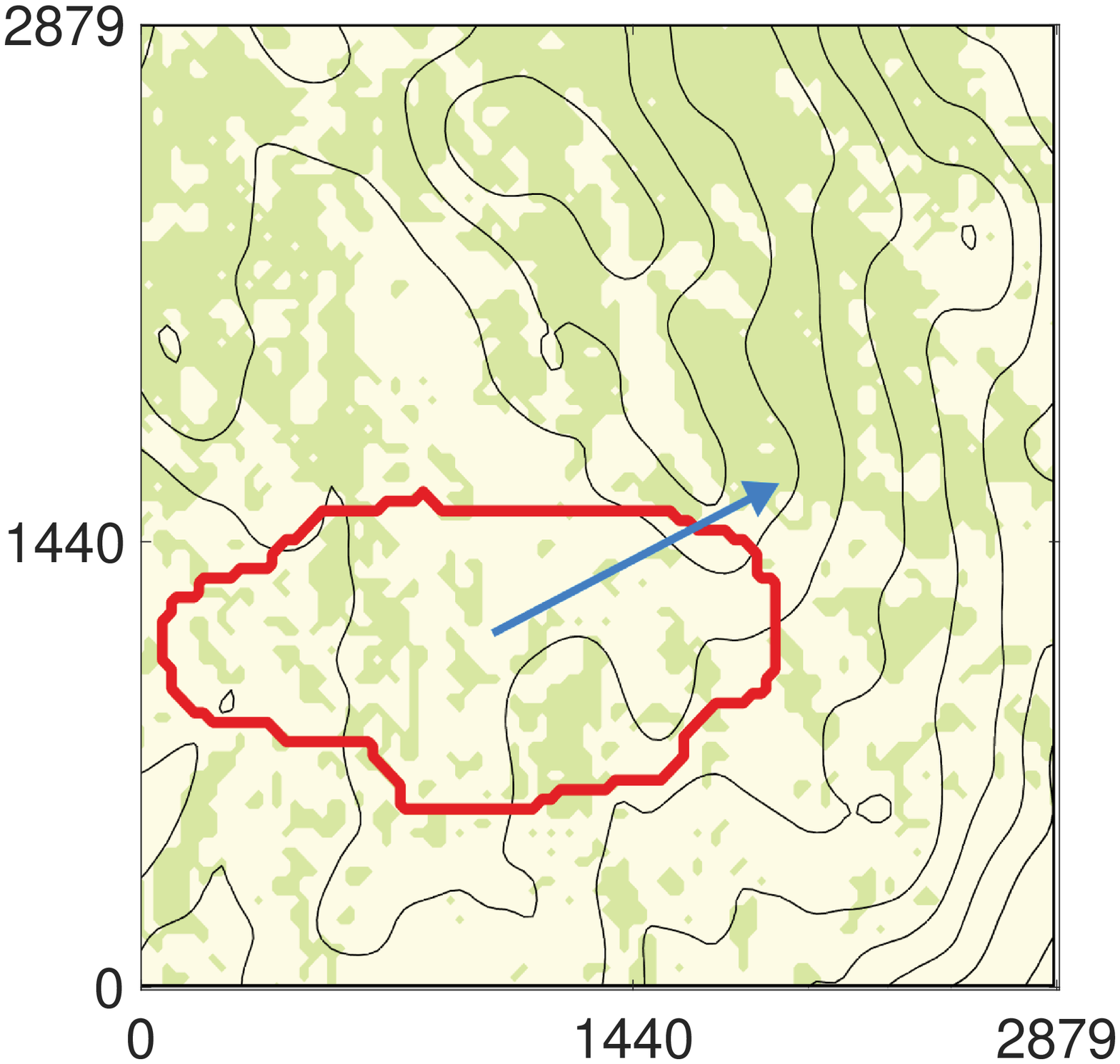}}\label{Troy_Initial.ps}%
\caption{Troy case study: 3D elevation map of the region where Troy fire propagated (a) and initial condition at 14:05 (b). The initial fire front is in red, the wind direction is in blue, the elevation level curves are in black, fuel A is in light-yellow, and fuel B is in light-green.}
\label{fig:troy}
\end{figure*}

\begin{figure*}[tb]
\centering
\subcaptionbox{}{\includegraphics[width=4.5cm]{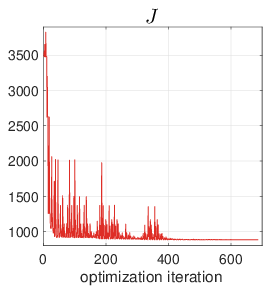}}%
\hfill
\subcaptionbox{}{\includegraphics[width=12.0cm]{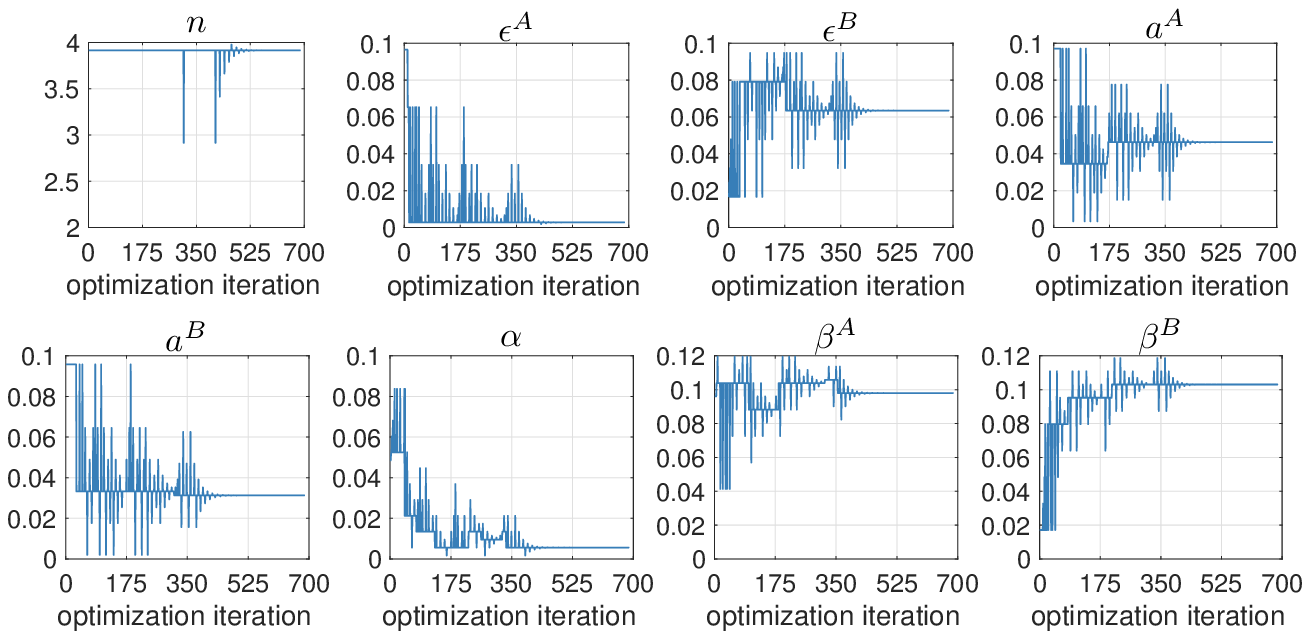}}%
\caption{Troy case study: trends of the cost function (a) and of the parameters (b) during the estimation procedure between 14:05 and 15:01.}
\label{fig:troyTrends}
\end{figure*}

\begin{figure*}[tb]
\centering
\subcaptionbox{}{\includegraphics[width=6.0cm]{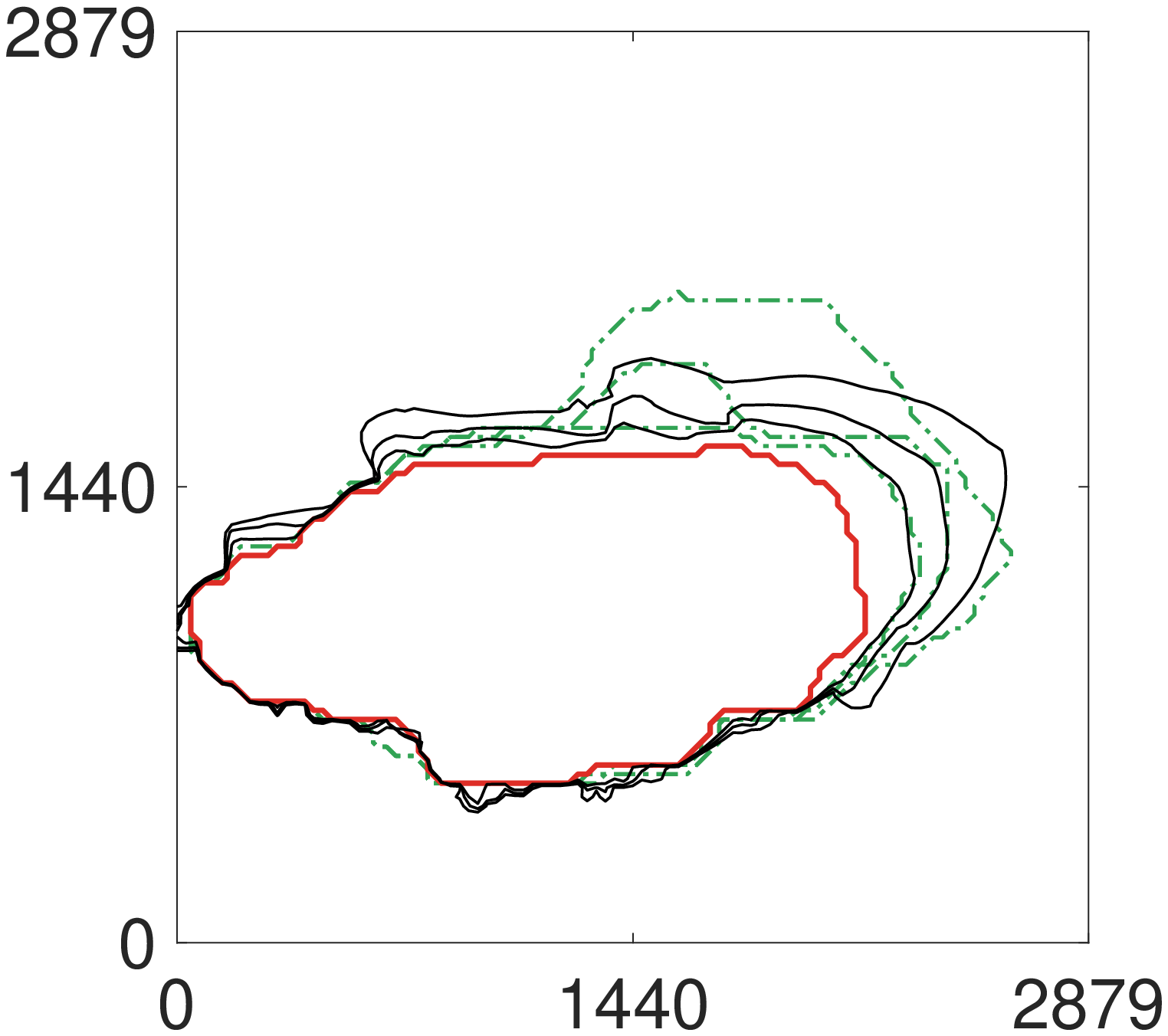}}
\subcaptionbox{}{\includegraphics[width=6.0cm]{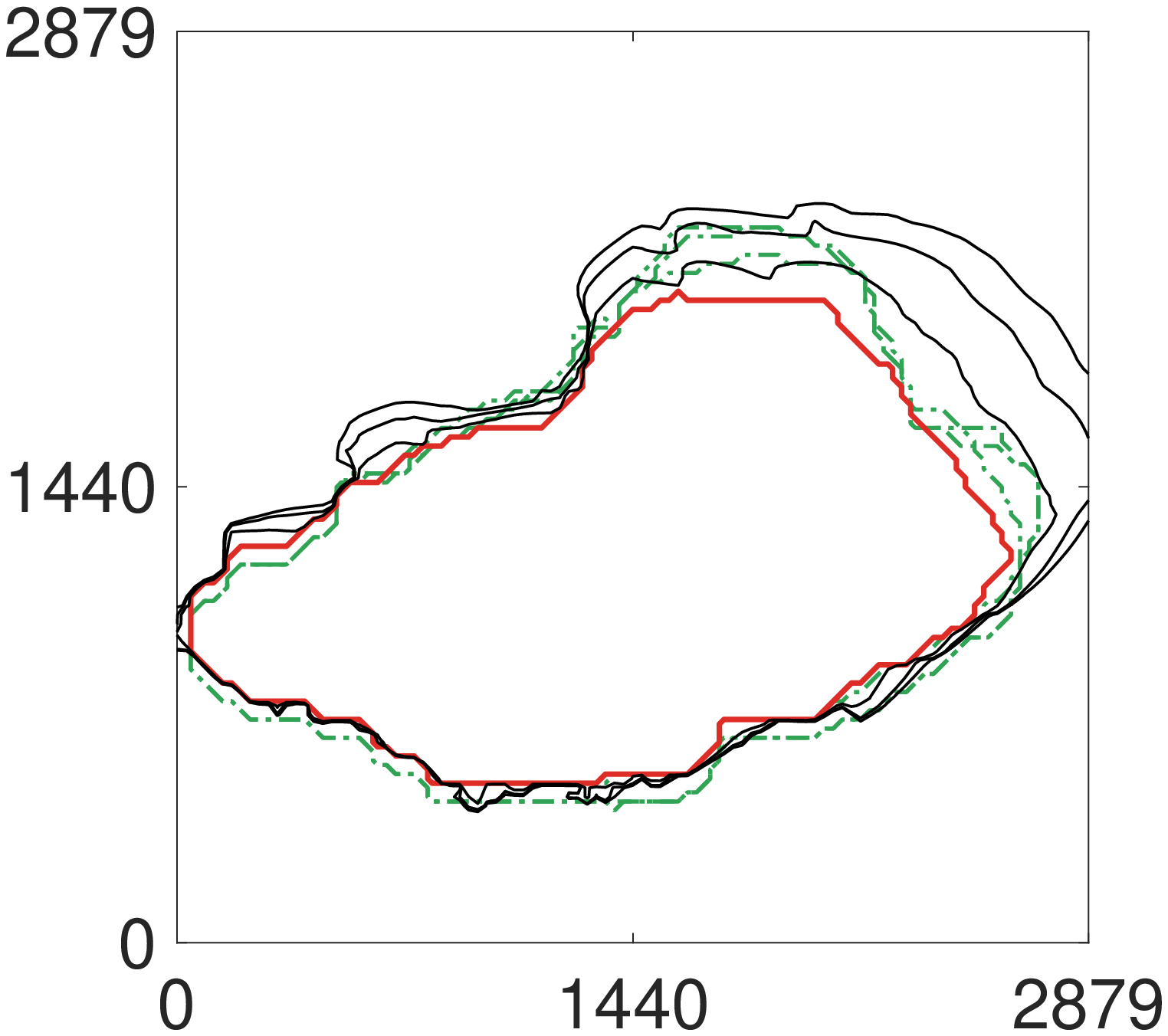}}
\caption{Comparison between simulation and real fronts of the Troy fire with $\v p=\vh p_1$. Fire evolution is compared at 14:05, 14:25, 14:39, and 15:01 in (a) and at 15:01, 15:21, 15:42, and 15:55 in (b). In both cases, the initial fronts at 14:05 or 15:01 are in red, simulated fronts are in black, and measured fronts are in dotted green.}
\label{fig:troyRealSim1}
\end{figure*}

\begin{figure*}[tb]
\centering
\subcaptionbox{}{\includegraphics[width=4.5cm]{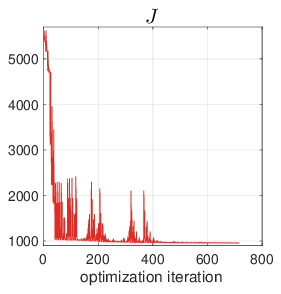}}%
\hfill
\subcaptionbox{}{\includegraphics[width=12.0cm]{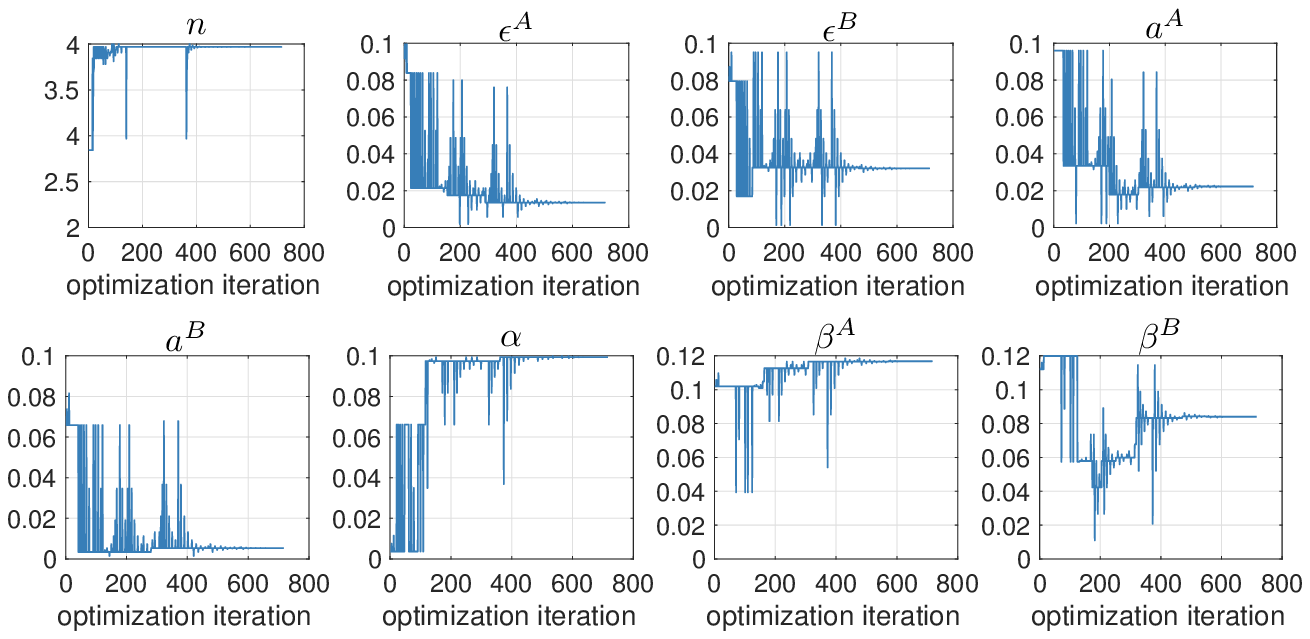}}%
\caption{Troy case study: trends of the cost function (a) and of the parameters (b) during the estimation procedure between 15:01 and 15:55.}
\label{fig:troyTrends2}
\end{figure*}

\begin{figure*}[tb]
\centering
\includegraphics[width=6.0cm]{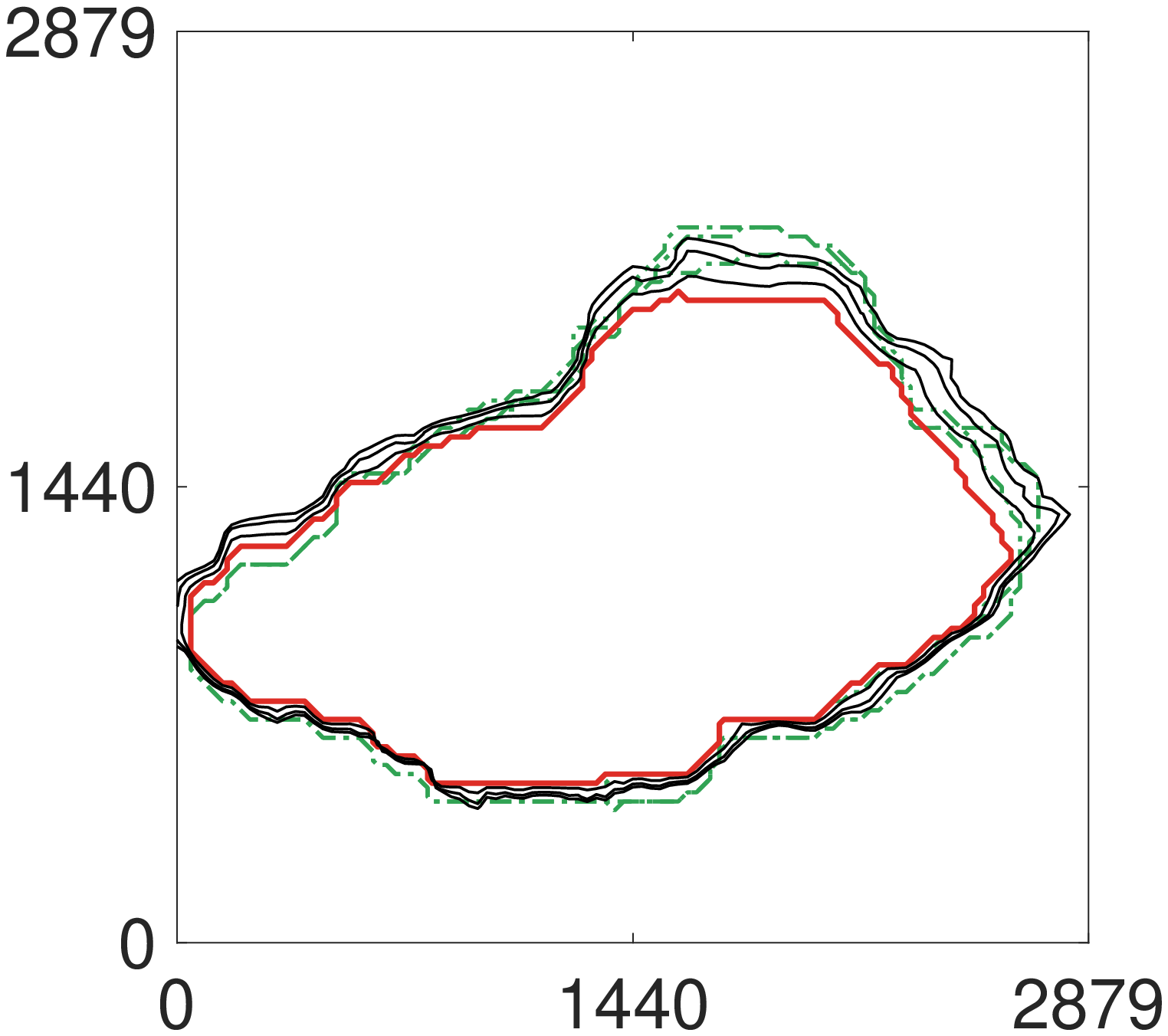}
\caption{Comparison between simulation and real images of the Troy fire with $\v p=\vh p_2$. Fire evolution is compared at 15:01, 15:21, 15:42, and 15:55. The initial front at 15:01 is in red, simulated fronts are in black, and measured fronts are in dotted green.}
\label{fig:troyRealSim2}
\end{figure*}

After testing the model on four known scenarios (see Section \ref{sec:sim:testing}) and verifying the effectiveness of the proposed parameter estimation approach in two simulated case studies (see Sections \ref{sec:sim:valley} and \ref{sec:sim:hill}), the proposed method is applied to a real case study. Toward this end, we consider the 2002 Troy fire. This fire started in Laguna Mountains in San Diego County, California, on June 19, 2002 and burnt about $1188$ acres of land \cite{Kim2012Thermal-image-basedEquation}. The evolution of this fire was recorded by Pacific Southwest Research Stations through measurements of the emitted thermal-infrared light from an aerial view, between 13:30 and 15:55 \cite{PacificSouthwestResearchStation2007TheScience}. From the observation of the aerial images (sketched in Fig. \ref{fig:troySatellite}), it is possible to notice a change of behavior occurring around 14:00. In fact, the wildfire front spreads mainly towards East and South-East between 13:30 and 13:46, while it deviates towards North-East after 14:05.

We start considering the fire fronts at 14:05, 14:25, 14:39, and 15:01 for a first round of test. Thus, the fire images captured by satellites are used to perform parameter identification, setting the front at 14:05 as initial condition. A squared area of $2879 \times 2879$ m$^2$ is taken as spatial domain, i.e., we choose $\Omega=[0,2879]\times[0,2879]$. Such a region is discretized with a regular grid made up of $101 \times 101$ cells. Local map and elevation data are obtained from Google Earth software and shown in Fig. \ref{fig:troy}(a). Vegetation data are collected from the US Forest Service LANDFIRE website \cite{U.S.DepartmentofAgricultureForestServiceandU.S.DepartmentoftheInteriorLANDFIRETools} and then grouped together in two main fuel categories, named A and B as in the previous two case studies. Fuel A includes short grass, timber grass, chaparral, and brush, while fuel B is composed of short needle litter, long needle or hardwood litter, timber litter, and understory. In order to apply the method proposed in this paper, it is necessary to know wind speed and direction during the whole fire evolution. Based on the observations of \cite{Lo2012AMethod} and on how the Troy fire evolved during the considered period, we fix wind velocity to a constant value equal to 30 mph, while the wind direction is supposed to be constant and equal to $63^\circ$ (being $0^\circ$ the North direction and $90^\circ$ the East direction). Fig. \ref{fig:troy}(b) shows the fire front at $t_0=\text{14:05}$, which, as said, is the initial condition of this test, together with elevation level curves, fuel distribution, and wind direction.

As previously pointed out, we first consider the fire fronts between time instants  $t_0=\text{14:05}$ and $t_f=\text{15:01}$ to perform parameter estimation. The objective function $J$ is evaluated at 14:25, 14:39 and 15:01, according to \eqref{eq:trackingCostH}. The initial value of the parameter vector $\v p$ is randomly chosen between $\v p_{\rm min} = (2, 0.001, 0.001, 0.001, 0.001, 0.0001, 0.0001, 0.0001)$, and $\v p_{\rm max}$ = (4, 0.1, 0.1, 0.1, 0.1, 0.1, 0.12, 0.12). The trends of the components of $\v p$ and of the cost function $J$ during the optimization procedure are shown in Fig. \ref{fig:troyTrends}. Optimization was performed in $13213$ s.  

Likewise in the valley and hill case studies, after a transient behavior, the cost and the components of the vector $\v p$ converge to stationary values, which represent the optimal estimates. The final estimate of the parameter vector $\v p$ obtained by the optimization algorithm is $\vh p_1= $ (3.91, 0.0028, 0.0635, 0.0463, 0.0313, 0.0056, 0.0980, 0.1030), while the optimal cost is equal to 882. A simulation of fire propagation with $\v p=\vh p_1$ from $t_0=\text{14:05}$ to $t_f=\text{15:01}$ provides the fire front snapshots depicted in Fig. \ref{fig:troyRealSim1}(a). The ratio $r$, computed as defined in \eqref{eq:r}, is equal to $0.0865$, while it is not possible to compute the error $e$ defined in \eqref{eq:e} since the value of the true parameters $\v p^*$ is not available.

To better evaluate the accuracy of parameter estimation, we introduce three additional performance indexes taken from the reference literature (see \cite{filippi2014evaluation}), i.e., the Sørensen similarity index (denoted by SSI), the Jaccard similarity coefficient (JSC, for short), and the Kappa statistics (named KS). They are defined as follows:
\begin{equation}
\text{SSI}(t,\v p) := \dfrac{2 |\Gamma^{\rm meas}(t) \cap \Gamma(t,\v p)|} {|\Gamma^{\rm meas}(t)|+|\Gamma(t,\v p)|},
\label{eq:sorensen}
\end{equation}
\begin{equation}
\text{JSC}(t,\v p) := \dfrac{|\Gamma^{\rm meas}(t) \cap \Gamma(t,\v p)|} {|\Gamma^{\rm meas}(t) \cup \Gamma(t,\v p)|},
\label{eq:jaccard}
\end{equation}
\begin{equation}
\text{KS}(t,\v p) := \dfrac{P_a(t,\v p) - P_e(t,\v p)} {1 - P_e(t,\v p)},
\label{eq:kappa}
\end{equation}
where $|\Gamma(t,\v p)|$ denotes the area of the region within the front $\Gamma(t,\v p)$, as well as
\begin{displaymath}
P_a(t,\v p) := \frac{|\Gamma^{\rm meas}(t) \cap \Gamma(t,\v p)|}{|\Omega|} + \frac{|\Omega \setminus (\Gamma^{\rm meas}(t) \cup \Gamma(t,\v p))|}{|\Omega|},
\end{displaymath}
\begin{displaymath}
P_e(t,\v p) := \frac{|\Gamma^{\rm meas}(t)| | \Gamma(t,\v p)|}{|\Omega|^2} + \frac{|\Omega \setminus \Gamma^{\rm meas}(t)||\Omega \setminus \Gamma(t,\v p)|} {|\Omega|^2}.
\end{displaymath}

Table \ref{tab:Troy1} contains the values of the performance indexes \eqref{eq:sorensen}, \eqref{eq:jaccard}, and \eqref{eq:kappa} computed on the basis of real images and results of simulations at 14:25, 14:39 and 15:01, with $\v p =\vh p_1$. Referring again to \cite{filippi2014evaluation}, we point out that the obtained values are large, thus proving the effectiveness of the proposed approach.

\begin{table}[t]
\centering
\caption{Values of performance indexes \eqref{eq:sorensen}, \eqref{eq:jaccard}, and \eqref{eq:kappa} for the Troy case study obtained between 14:05 and 15:55 with $\v p =\vh p_1$.}
\label{tab:Troy1}
\begin{tabular}{c|ccc|ccc} 
\toprule
& $t\!=\!\text{14:25}$ & $t\!=\!\text{14:39}$ & $t\!=\!\text{15:01}$ & $t\!=\!\text{15:21}$ & $t\!=\!\text{15:42}$ & $t\!=\!\text{15:55}$ \\
\midrule
$\text{SSI}(t,\v p)$ & 0.9416 & 0.9082 & 0.8439 & 0.9221 & 0.8190& 0.7858 \\
$\text{JSC}(t,\v p)$ & 0.9613 & 0.9362 & 0.8806 & 0.9396 & 0.8442 & 0.8078 \\
$\text{KS}(t,\v p)$ & 0.9699 & 0.9519 & 0.9153 & 0.9595 & 0.9005 & 0.8801 \\
\bottomrule
\end{tabular}
\end{table}

\begin{table}[t]
\centering
\caption{Values of performance indexes \eqref{eq:sorensen}, \eqref{eq:jaccard}, and \eqref{eq:kappa} for the Troy case study obtained between 15:01 and 15:55, setting $\v p =\vh p_2$.}
\label{tab:Troy2}
\begin{tabular}{c|ccc} 
\toprule
& $t\!=\!\text{15:21}$ & $t\!=\!\text{15:42}$ & $t\!=\!\text{15:55}$ \\
\midrule
$\text{SSI}(t,\v p)$ & 0.9424 & 0.8957& 0.8987 \\
$\text{JSC}(t,\v p)$ & 0.9564 & 0.9165 & 0.9178 \\
$\text{KS}(t,\v p)$ & 0.9703 & 0.9450 & 0.9466 \\
\bottomrule
\end{tabular}
\end{table}

After estimating the parameters $\vh p_1$, we perform a further simulation run with the identified parameters $\vh p_1$ in order to evaluate the capability to forecast the evolution of the fire front. In more detail, we let $\v p=\vh p_1$, as well as we set the initial fire front shape equal to the real front at 15:01 and fixed the final time at 15:55. The simulation results are then compared with the real fire front shapes at 15:21, 15:42 and 15:55. The results are reported in Fig. \ref{fig:troyRealSim1}(b). It is possible to observe an over-estimation of the fire propagation towards North-East, which brings a small part of the fire front outside of the spatial domain, and a good approximation of the fire evolution in the other directions. The performance indexes \eqref{eq:sorensen}, \eqref{eq:jaccard}, and \eqref{eq:kappa} of this simulation run are shown again in Table \ref{tab:Troy1}. They are still large, even if we may observe a reduction of their values, especially at 15:55.

The last result suggests that the accuracy of the forecasts decreases if we consider a constant set of parameters over a large time span, possibly due to changes in the wildfire regime. Thus, after postulating that a periodical re-run of the parameter estimation procedure may increase the effectiveness of estimation of the fire fronts, we performed a second round of tests by estimating parameters using measurements of the fire front between $t_0=\text{15:01}$ and $t_f=\text{15:55}$. The objective function $J$ is evaluated at 15:21, 15:42 and 15:55, according to \eqref{eq:trackingCostH}. The trends of the components of $\v p$ and of the cost function $J$ during the optimization procedure are shown in Fig. \ref{fig:troyTrends2}. Optimization was performed in $7655$ s.  
The final estimate of the parameter vector $\v p$ obtained by the optimization algorithm is $\vh p_2=$ (3.97, 0.0135, 0.0321, 0.0223, 0.0054, 0.0994, 0.1168, 0.0840), while the optimal cost is equal to 958. The value for the ratio $r$ in \eqref{eq:r} is equal to $0.0939$. The verify the effectiveness of the estimation, we let $\v p=\vh p_2$ and simulate the wildfire evolution from 15:01 to 15:55.  Fig. \ref{fig:troyRealSim2} shows the obtained fire front shapes. If we compare these results with those depicted in Fig. \ref{fig:troyRealSim1}(b), we notice that the fire fronts obtained with $\vh p_2$ are significantly more similar to the real fronts, with no over-estimation of the fire evolution towards North-East. The larger values of the performance indexes \eqref{eq:sorensen}, \eqref{eq:jaccard}, and \eqref{eq:kappa} in Table \ref{tab:Troy2} confirm this result from a quantitative point of view. Therefore, we conclude that a periodical re-run of the parameter estimation can improve the simulation accuracy, taking into account possible changes of the fire regime.

\section{Conclusions}
\label{sec:concl}

In this paper, we have addressed the problem of modeling fire propagation using level set methods to describe the front evolution over space and time and an empirical model for the rate of spread. The model is composed of several parameters depending on vegetation and characteristics of  terrain, which have been estimated by means of a procedure based on the minimization of a least-squares cost accounting for the difference between measured and predicted fronts. The proposed combination of modeling and identification has been tested with both simulated data and real measurements of propagating fronts. In particular, we have considered data of the 2002 Troy fire, and we have obtained optimal estimates of the parameters of the propagating model that have turned out to be effective in predicting the future evolution of the fire front.

The successful results we have obtained suggest several future directions of investigation. The first one is the use of other optimization techniques and the adaptation over time of model parameters, which may be well-suited to predicting fire propagation for the purpose of decision support to allocate the firefighting resources and contain the fire. Moreover, we will extend our approach to more complex models of the rate of spread than the one proposed in \cite{Mallet2009ModelingMethods, Lo2012AMethod}. Another important subject of future effort will be the validation of the considered fire propagation models with estimated parameters. Thus, according to \cite{Morvan2019}, we will test the identified models without variations of some parameters against several real case studies in order to better detect case-dependent and case-independent parameters as well as quantify the quality of the overall estimation.

%\section*{References}

\bibliography{references}

\end{document}